\documentstyle[12pt,amsfonts]{article}
\newfam\msafam
\newfam\msbfam
\newfam\euffam

\newcommand{\A}{{\Bbb A}}

\newcommand{\C}{{\Bbb C}}

\renewcommand{\H}{{\Bbb H}}

\newcommand{\N}{{\Bbb N}}

\renewcommand{\P}{{\Bbb P}\hspace{.05em}}

\newcommand{\Z}{{\Bbb Z}}

\newtheorem {Theorem}{Theorem}[section]
\newtheorem {Proposition}[Theorem]{Proposition}
\newtheorem {Corollary}[Theorem]{Corollary}
\newtheorem {Definition}[Theorem]{Definition}

\newtheorem {Lemma}[Theorem]{Lemma}
\newtheorem {Remark}[Theorem]{Remark}
\newtheorem {Notation}[Theorem]{Notation}
\newtheorem {Claim}[Theorem]{Claim}

\newcommand{\prf}{{\it Proof:} }
\newcommand{\qed}{\hspace*{\fill}\hbox{$\Box$}\\ \ \\}

\newcommand{\XX}{\hbox{$S \times {\P}^n$}}   

\newcommand{\smX}{\hbox{\small $S \times \P^n$}}   
\newcommand{\smXS}{\hbox{\small $S \times  \P^n /S$ }} 
   
\newcommand{\smXSdos}{\hbox{\small $S \times \P^2/S$}}  
\newcommand{\smholS}{\hbox{\small ${\cal O}_S$}}

\newcommand{\arrA}{\hbox{${\cal A}$}}       
\newcommand{\coneA}{{\bf c}\arrA}        
\newcommand{\hol}{\hbox{$\cal O$}}     
\newcommand{\flecha}{\longrightarrow} 

\newcommand {\pn}{{\P}^n}
\newcommand {\pndual}{{\P}^{n^\vee}}   
\newcommand {\join}{\vee}

\newcommand{\calL}{\hbox{${\cal L}$}}

\newcommand{\calC}{\hbox{${\cal C}$}}
\newcommand{\calE}{\hbox{${\cal E}$}}
\newcommand{\calJ}{\hbox{${\cal J}$}}

\newcommand{\calM}{\hbox{${\cal M}$}}           
\newcommand{\bfF}{\hbox{${\bf F}$}}           

\begin{document}
\begin{center}{{\Large 
Monodromy of hypergeometric functions\\
arising from arrangements of hyperplanes}\\ \ \\
Herbert Kanarek\\
Universit\"at GH-Essen\\
FB6 Mathematik\\
45117 Essen, Germany}
\end{center}


\vskip1cm

\centerline{\large\bf  Abstract}
\vskip.5cm

Given an arrangement of hyperplanes in $\P^n$, possibly with
non-normal crossings, we give a vanishing lemma for the 
cohomology of the sheaf of $q$-forms with logarithmic poles
along our arrangement. We give a basis for the ideal $\cal J$ of relations
for the Orlik-Solomon's algebra.  Under certain genericity conditions
it was shown by H.~Esnault, V.~Schechtman and E.~Viehweg that the
cohomology of a local system is given by the Aomoto complex.
We generalize this result to a deformation of local systems obtained
via a deformation of our arrangement.  We calculate the
Gau\ss-Manin connection for this case.  We give a basis for the
Gau\ss-Manin bundle for which, with help of the basis 
for $\cal J$, we give then a method to calculate a representation 
of this connection.  From here, with the results of K-T.~Chen or 
P.~Deligne, one can calculate the monodromy representation.  This
gives a generalization of the hypergeometric functions.


\section*{Introduction}
\label{introduction articulo}

Let $\{H_i\}_{i\in I}$ be a collection of different hyperplanes in $\pn$,
let $U=\pn\setminus \{H_i\}_{i\in I}$ and
let $X$ be a smooth compactification of $U$ in such 
a way that the divisor $D=X\setminus U$ has normal crossings.
Let $\omega$ be a holomorphic ${\it gl}(n,U)$-valued $1$-form with
logarithmic poles along $D$.  One can define 
a holomorphic connection $\nabla$
on $\hol_U\otimes \C^n$ as
\[
\nabla = df + f\omega.
\]
On $\Omega^p_U\otimes \C^r$ one defines $\nabla$ by 
$\nabla(\alpha\otimes v) = d\alpha \otimes v + (-1)^p\alpha \nabla(v)$.
If $d\omega - \omega\wedge\omega = 0$ then 
$\nabla^2 = 0$ and one says that the connection is integrable.  
Let $V$ be the rank $r$ system of flat sections
of $(\hol_U^r,\nabla)$, i.e. $V$ is the kernel of $\nabla$ which locally
analytically is isomorphic to $\C^r$.  

In this paper we attempt to describe the monodromy of hypergeometric
functions arising as solutions of ordinary differential equations
with regular singularities along an arrangement of hyperplanes 
in the projective space.
As it is well known, from Riemann's integral representation formula one can 
express the hypergeometric system of differential equations as
a direct image (as a variation of cohomology) of a rank one system,
see [S, Theorem 2] and [M, Proposition in p.~373]. We
consider a deformation of rank one systems on the complemet of
a configuration in $\pn$ and obtain this generalization of 
the hypergeometric function. We use the following results 
on rank one systems.

Consider a collection of different hyperplanes  in $\pn$ as before.
Let $\omega$ be a global holomorphic differential form over
$U=\pn\setminus \{H_i\}_{i\in I}$ such that it has logarithmic 
poles along a divisor $D$.  As before, this form induces an integral 
connection $\nabla$ on $\hol_X$ which as flat section gives a 
rank 1 local system. Deligne proved in [D1] that if one takes $\omega$ in 
such a way that it has no positive integers as residues then the 
cohomology of the local system of flat sections of $\nabla$ is given 
by the cohomology of the de Rham complex induced by $\nabla$. Moreover, 
if we suppose that $\omega$ has no integers as residues, Esnault, 
Schechtman and Viehweg showed in [ESV] that the cohomology of the 
local system is given by the Aomoto complex which is the subcomplex 
of global sections of the de Rham complex.
Under the same genericity conditions, Esnault and Viehweg 
proved in [EV] that the cohomology of $V$ is concentrated in the
$n$-th term. Actually all these results are not only given for
higher rank systems.

To describe the monodromy we use the previous results and consider 
deformations of local systems. We take a connection over a rank one 
bundle with logarithmic
poles along $\arrA$, as before, and make deformation 
of this connection in such a way that all the residues remain 
constant.  For this we consider a topologically trivial deformation of
the arrangement $\arrA = \sum_{i\in I} H_i \subset \pn$.  This 
deformation is given by moving one hyperplane and leaving 
the other ones fixed in such a way that we don't get new non-normal 
crossings.  Let $H_{i_0}$ with $i_0\in I$ be such hyperplane.  We then
have a family of arrangements parametrized by the complement of 
the discriminant of the arrangement 
$\arrA^\prime = \cup_{i\in I\setminus {i_0}}H_i$ which is 
given by our original arrangement without the hyperplane that we 
want to move.  We obtain a relative connection on our family of 
arrangements.  We show that one can extend this connection to an 
absolute connection $\overline \nabla$.  In this way, and it is the 
deciding point, one can construct  the Gau\ss-Manin connection.  
This connection is defined on the Gau\ss-Manin bundle which is defined 
as the relative de~Rham cohomology sheaves and has as flat section
the direct images of the absolute local system on the family. 
Theorem \ref{teo:esv-ext NNC} is a generalization of the results 
obtained by Esnault, Schechtman, Viehweg concerning the cohomology 
of the family.  They imply that the Gau\ss-Manin bundle is generated 
by global sections.  One can then, in an standard way, calculate 
a representation of the Gau\ss-Manin connection.  
In the general case, to calculate the monodromy, one can use 
the ideas of K-T.~Chen of iterated integrals to calculate the solutions
of our system of differential equations but in Example I we 
rather prefer to use other results on differential equations, 
see [D1, II.5.6].

In the first section we give a vanishing lemma for the cohomology 
of differential forms with logarithmic poles along
our arrangement.  We will use this lemma to show that the
cohomology of the local system is given by the Aomoto complex.
The second section is devoted to give a description
of the Gau\ss-Manin connection.  We consider a 
family of arrangements parametrized by the complement of the 
discriminant of an arrangement in $\pn$.  Given 
a relative connection we construct the Gau\ss-Manin connection 
with logarithmic poles along the divisor obtained from a 
desingularization of the divisor $Discr(\arrA) \subset \pndual$.

In the third section contains several results concerning some aspects of 
the combinatorics for the theory of hyperplane
arrangements.  We present two important results. The first one due to
Bj\"orner, see [Bj], where he gives a base for the cohomology 
with constant coefficients on the complement of an arrangement. By 
a result of Deligne [D2] this cohomology is equal to the cohomology of
the sheaf of $q$-differential forms with poles along an arrangement.
The second result is a basis for the ideal of relations
$\calJ$ for the Orlik-Solomon's Algebra.  

In section 4 we have the principal results. We show in Theorem 
\ref{teo:esv-ext NNC} that for non normal crossing case the 
cohomology of the relative local system is also given by the Aomoto 
relative complex. From here one can give a representation of the
Gau\ss-Manin connection. In section 5 and 6 we give some examples.
We construct the Gau\ss-Manin 
connection for the following arrangement:
\[
\begin{array}{c}
\unitlength1cm
\begin{picture}(7,5)
\put (1.7,0.7){\line(1,1){2.7}}
\put (2.3,1.9){$\small{H_5}$}
\put (1.4,3.3){\line(2,-1){5}}
\put (2.3,2.9){$\small{H_4}$}
\put (1.8,5.4){\line(1,-2){2.4}}
\put (3.2,1.5) {$\small{H_3}$}
\put (1.5,1){\line(1,0){5}}        
\put (5, 0.5) {$H_2$}
\put (2,0.4){\line(0,1){5}}        
\put (1, 2 ) {$\tiny{H_1}$}
\put (1.7,5.3){\line(1,-1){4.6}}
\put (4.7,2.5) {$H_0$}
\put (3, 0){$Fig.~a$} 
\end{picture}
\end{array}
\]
This example is of particular interest because as discriminant one
obtains the ``Ceva'' configuration which has been deeply studied in
[BHH] for the construction of Ball quotient surfaces. The second
example is an example in which the arrangement chosen has non-normal 
crossings.

I am very grateful to H.~Esnault for her help and encouragement
during the research and preparation of this paper.  I would also
like to thank E.~Viehweg for his help and useful suggestions,
in particular during the preparation of
the first section.  I wish express my gratitude to R.~Hain, 
S.~M\"uller-Stach and to V.~Welker for the useful discussions 
and suggestions, in particular Proposition \ref{prop:rel} 
was obtained together with V.~Welker.


\section{Some Vanishing Theorems}
\label{sec. Vanishing Theorems}

Let $ \{H_i\}_{i\in I}$  be a family of distinct hyperplanes in 
$\P^n$,  $H = \sum_{i \in I} H_i$  the associated effective 
divisor and $U = \pn \setminus H$ be the complementary affine open set.
We have the following
definition given in [ESV, Definition~({\bf Bad})].\\

\begin{Definition} 
\label{def. bad}
a) Given a linear subspace $L \subset \P^n$ let
$$I_L = \{i \in I ~|~ L \subset H_i \}.$$
b) We define the set
\[
\begin{array}{c}
{\cal L}_j(H) = \{ L \subset \P^n~{\rm linear }~|~ 
         dim L = j~{\rm and }\\
         L = \cap _{i \in I_L \backslash \{ i_0 \} } H_i~
         {\rm for~every }~i_0 \in I_L \}
\end{array}
\]
for $0 \leq j \leq n-2.$
Let
$$\calL(H) = \cup _{j=0}^{n-2} {\cal L}_j(H).$$
\end{Definition}
The loci where $H$ has non-normal crossings are
exactly the linear subspaces contained in ${\cal L}(H)$.  When
there is no possible confusion about the divisors we will write 
only $\calL$.  
Let $X$ be the variety obtained by considering successive blow ups 
along the elements of ${\cal L}$ in the following way.

Let $\pi ^{(r)} = \tau_1 \circ ...\circ \tau_r$

\begin{equation}
X_r \stackrel{\tau_r}{\longrightarrow} X_{r-1} 
\stackrel{\tau_{r-1}}{\longrightarrow}\ldots\flecha
X_1\stackrel{\tau_1}{\longrightarrow} \pn
\end{equation}
where $\tau _i$ is the blow up of $X_{i-1}$ along the proper 
transform $T_{i-1}$ under 
$\pi^{(i-1)}$ of the elements of ${\cal L}_{i-1}$.  
As shown in [ESV, Claim], $T_{i-1}$ is the disjoint union
of closed nonsingular subschemes.
Writing $X_0 = \P^n$ and $X = X_{n-1}$
let us set for $r\leq n-1$
\[
\pi _{(r)} = \tau_{r+1} \circ ... \circ \tau_{n-1} : X = X_{n-1} 
\longrightarrow X_r
\]
and
\[
\pi = \pi_{(0)}  = \pi^{(n-1)}.
\]
The variety $X$ is nonsingular and $\pi$ will be called a {\it standard} 
resolution of $H$.

\begin{Lemma}
\label{lem:blup} 
Let $I' \subset I$ and let us consider the divisor $H' = 
\sum_{i \in I'} H_i$.  
Let \newline
$\pi \prime : X' \longrightarrow \P^n $ be the standard
resolution of $H'$. Then there exists a morphism 
$\gamma : X \longrightarrow X' $ such that
\[
\begin{array}{lcr}
       X          &  \stackrel{\gamma}{\longrightarrow} &          X'\\
   \pi \searrow   &                                     &\swarrow \pi'\\
                  &              \P^n                   &       
\end{array}
\]
commutes.\\
\end{Lemma}

\prf
Let $\calL'_j$ be the bad strata of dimension $j$ of
$H'$ and let $\calL' = \cup_{j=0}^{n-2}\calL'_j$.
Note that ${\calL'_j \subset \calL_j}$. Let 
$\tau^\prime_j : X'_j \flecha X'_{j-1}$ be the $j$-th blow up 
of $H'$ and let $\pi'^{(j)}= \tau'_1\circ\dots\circ\tau'_{j}$. 
The proof follows by induction over the dimension of the bad loci.
\qed

We have two special cases in which we would like to apply the previous
lemma.\\

a) Let $L \in {\cal L}_j$ and $I' = I_L$, i.e.,  
     $$I' = I_L = \{ i \in I | L \subset H_i \}$$
so we have the commutative diagram
\[
\begin{array}{lcr}
X              & \stackrel{\gamma}{\longrightarrow}  &      X_L\\
 \pi  \searrow &                                     & \swarrow \pi _L\\
	       &              \pn                    &           
\end{array}
\]
where we write $X_L = X'$ and $\pi_L=\pi'$. In this case the
exceptional divisor  
$E_L = \pi^*(L)$ is a $\P^{n-j-1}$ bundle over $L$ for which
on every fiber, the proper transform of $H_L = \cup_{i\in I_L}H_i$ 
gives a configuration 
$\Gamma = \sum_{i \in I_L} \Gamma_i \subset \P^{n-j-1}$.

On the other hand, for $L\in {\cal L}_j$ we have a 
configuration in $L \cong \P^j$ given by
\[ 
\sum_{i \in I\backslash I_L}H_i\mid _L .
\]

b) For any $i_0 \in I$ let $I' = I \backslash \{i_0\}$ then we have 
\[
\begin{array}{lcr}
X             & \stackrel{\gamma}{\longrightarrow} &       X'\\
 \pi \searrow &                                    & \swarrow \pi' \\
              &                 \pn                &          
\end{array}.
\]

\begin{Claim}
\label{claim:biggerresolutions}
Let $H= \sum_{i \in I} H_i$ as before, $\pi:X \flecha \pn$ the
standard resolution and $D=\pi^*(H)$. Let $\varpi: Z\flecha \pn$ 
be any other resolution for which there 
exists a morphism $\gamma : Z \flecha X$ with $\pi\circ\gamma = \varpi$.
If $B= \gamma^*(D)$, then for $\nu \geq 0$ and for all 
$p, q \geq 0$ we have
\[
H^p(Z,\Omega_Z^q(\log B)\otimes\gamma^*\pi^*(\hol(\nu))) =
\]
\[
  H^p(X,\gamma_*\Omega_{Z}^q(\log B)\otimes\pi^*(\hol(\nu)))=
   H^p(X,\Omega_X^q(\log D)\otimes\pi^*(\hol(\nu))).
\]
\end{Claim}
\prf By [EV, Lemma 3.22] we know that for $p>0$
$R^p\gamma_*\Omega_Z^q(\log B) = 0$ and that 
$\gamma_*\Omega^q_Z(\log B) = \Omega^q_X(\log D)$.
Applying the projection formula we have 
\[
R^p\gamma_*(\Omega_Z^q(\log B)\otimes\gamma^*\pi^*\hol(\nu))
=R^p\gamma_*\Omega_Z^q(\log B)\otimes\pi^*\hol(\nu) = 0 
\hbox{\ \ \  for\ \ } p>0.
\]
We can now apply the Leray spectral sequence to obtain the first 
equality. The second equality is clear, since both sheaves are
isomorphic.  
\qed

Considering $H = \sum_{i \in I} H_i$ as before and 
$\pi : X \longrightarrow \pn$
the standard resolution along $H$, we have the following lemma.

\begin{Lemma} 
\label{lemmaVanishing}
Let $H=\sum_{i\in I}H_i$ be a non trivial configuration of 
hyperplanes in $\pn$, $\pi: X \flecha \pn$ a standard resolution and 
$D = \pi^*(H)$ be the reduced pull back divisor of $H$ then, 
for $p > 0$, $\nu \geq 0$, we have 
\[
H^p(X, \Omega_X^q(\log D) \otimes \pi^* \hol_{{\P}^n}(\nu )) = 0.
\]
\end{Lemma}

\prf  The proof will be by double induction over the number of 
hyperplanes and over $\nu$.
Let $H = H_1.$ We will consider first the case when $\nu =0.$

As for $k > 0$ we know that
\[
H^k({\A}^n, \C) = 0,
\]
by the degeneration of the Hodge to de~Rham spectral sequence,
see [D2, Corollary~3.2.13], we have that for $p+q>0$
\[
H^p(\P^n,\Omega^q(\log H_1)) = 0 .
\]
In this case one can argue by considering the long exact
sequence of cohomology associated to 
the short exact sequence
\begin{equation}
\label{eq:diff.form.seq}
0\flecha \Omega_{\P^n}^q\flecha \Omega_{\P^n}^q(\log H_1)
\otimes\hol_{\pn}\flecha \Omega^{q-1}_{\P^{n-1}}\flecha 0.
\end{equation}
The connecting morphism for the long exact sequence of 
cohomology obtained from (\ref{eq:diff.form.seq})
is an isomorphism.
As
\[
H^p(\P^n,\Omega^q_{\P^n})= \left\{ \begin{array}{c} 
0 \hbox{\ \ \ \ for \ \ } q \neq p \\
\C \hbox{\ \ \ \ for \ \ } q=p
\end{array} \right.
\]
we have the result.

For the case where $\nu > 0$ we tensor the sequence (\ref{eq:diff.form.seq})
by $\hol_{\pn}(\nu)$ to obtain the sequence
\begin{equation}
\label{eq:twistedseq}
0\flecha \Omega_{\P^n}^q(\nu)\flecha \Omega_{\P^n}^q(\log
H_1)\otimes\hol_{\pn}(\nu)\flecha 
\Omega^{q-1}_{\P^{n-1}}(\nu)\flecha 0.
\end{equation}
By Bott's formula, see [OSS, p. 8], we have 
that for $\nu>0$ and $p>0$
\[
H^p(\P^r,\Omega^q_{\P^r}(\nu)) = 0.
\]
Thus from $(\ref{eq:twistedseq})$ we have
\[
H^p(\pn,\Omega^q_{\pn}(\log H_1)\otimes\hol_{\pn}(\nu)) = 0
\hbox{\ \ \ \ for \ \ } p > 0
\]
which proves the lemma for the case when $\mid I \mid = 1$.

Let $I'\subset I$ be a proper subset. As induction hypothesis we
can assume that for $H'=\sum_{i\in I'}H_i$ 
\[
H^p(X',\Omega^q_{X'}(\log D')\otimes\pi'^*\hol_{\pn}(\nu)) = 0 
\]
where $\pi':X'\flecha \pn$ is the standard resolution of 
$H'=\sum_{i\in I'}H_i$ and $D'=\pi'^*(H')$.
Let $i_0 \in I$ fixed and $I'=I\setminus \{i_0\}$, 
by Lemma~$\ref{lem:blup}.(a)$ we have the
morphism 
\[
\gamma: X\flecha X'
\]
such that, for $\pi:X\flecha\pn$ the standard resolution
of $H$, the following diagram 
\[
\begin{array}{lcr}
\label{diagramainduction}
X             & \stackrel{\gamma}{\longrightarrow} &       X'\\
 \pi \searrow &                                    & \swarrow \pi' \\
              &                 \pn                &          
\end{array}
\]
is commutative.
Let $D= \pi^*(H)$ be the pull back of $H$ under $\pi$, 
$D' = {\pi'}^*(H')$ the pull back of $H^\prime$ under $\pi^\prime $
and $D'' = \gamma^* D^\prime$ the pull back of $D^\prime$ 
under $\gamma$ all taken to be reduced. Let
$D_{i_0} =D - D''$ which is equal to the proper transform of $H_{i_0}$ 
under $\pi$.  We have from [EV, 2.3.b] the following 
exact sequence
\begin{equation}
\label{suc.pullback-formasdif.1}
0 \flecha \Omega^q_X(\log\gamma^*(D')) \flecha \Omega^q_X(\log D) 
   \flecha \Omega^{q-1}_{D_{i_0}}(\log \gamma^*(D')\mid_{D_{i_0}}) 
   \flecha 0.
\end{equation} 
Tensoring with $\pi^*\hol(\nu)$ gives
\begin{equation}
\label{suc.pullbackformasdif}
\begin{array}{c}
0 \flecha \Omega^q_X(\log D'')\otimes\pi^*\hol(\nu) 
   \flecha
\\ \ \\
  \Omega^q_X(\log D)\otimes\pi^*\hol(\nu)
 \flecha \Omega^{q-1}_{D_{i_0}}(\log D''\mid_{D_{i_0}})
      \otimes\pi^*\hol(\nu)\flecha 0.
\end{array}
\end{equation} 
Applying Claim~\ref{claim:biggerresolutions}, by induction 
on the dimension 
\[
H^p(D_{i_0},\Omega^{q-1}_{D_{i_0}}(\log D''\mid_{D_{i_0}})\otimes
      (\pi\mid_{D_{i_0}})^*\hol(\nu))=0.
\]
In fact, $D_{i_0}$ is a resolution of a configuration in 
$\P^{n-1}\cong H_{i_0}$.  It is easy to see, that $D_{i_0}$
is the standard resolution for this configuration, but by 
Claim~\ref{claim:biggerresolutions}, this is not necessary
for the vanishing.

On the other hand, again by projection formula and
Claim~\ref{claim:biggerresolutions} we have
\begin{equation}
H^p(X, \Omega_X ^q(\log D'') \otimes \pi^*\hol_{\pn}(\nu)) =
H^p(X',\Omega^q_{X'}(\log D')\otimes 
\pi'^*\hol_{\pn}(\nu)). 
\end{equation}
By our induction hypothesis on the number of hyperplanes
both groups are zero, for $p>0$. 
From the long exact sequence of cohomology obtained 
from~(\ref{suc.pullbackformasdif}) we have
\[
H^p(X,\Omega^q_X(\log D)\otimes \hol_{\pn}(\nu)) = 0 
\hbox{\ \ \ \ for \ \ } p>0.
\]
\qed

As an interesting application of this lemma we have:

\begin{Corollary}
Let $H=\sum_{i\in I}H_i$ and $\pi:X \flecha \pn$ as 
in Lemma~\ref{lemmaVanishing}, for $p>0$ we have 
\begin{equation}
R^p\pi_*(\Omega^q_X(\log D)) = 0
\end{equation}
\end{Corollary}

\prf
As $\Omega_X^q(\log D)$ is coherent, by being locally free, and
as $\pi$ is proper, we have that 
$R^p\pi_*(\Omega^q_X(\log D)\otimes \pi^*\hol(\nu))$ is coherent,
see [Ha, Theorem III 8.8].
From Serre's Vanishing Theorem, there exist 
$\nu_0$ such that for every $\nu \geq \nu_0$ 
\begin{equation}
H^p(\pn, R^j\pi_*(\Omega^q_X(\log D)\otimes \pi^*\hol(\nu)))=0
\end{equation}
for $p>0$.
By the Leray spectral sequence we then have
\begin{equation}
\label{Leray on global sections}
H^0(\pn, R^j\pi_*(\Omega^q_X(\log D)\otimes\pi^*\hol(\nu)))=
H^j(X, \Omega^q_X(\log D)\otimes\pi^*\hol(\nu)) .
\end{equation}

On the other hand, as the sheaf $\hol(\nu)$ is ample, there exists 
$\nu_1$ such that for every
$\nu\geq \max\{\nu_0,\nu_1\}$ we have that 
$R^j\pi_*(\Omega^q_X(\log D))\otimes\hol(\nu)$ is generated
by global sections. By Lemma~\ref{lemmaVanishing} 
$H^p(X, \Omega^q_X(\log D)\otimes\hol(\nu))=0$. This implies, via
(\ref{Leray on global sections}), that 
$R^j\pi_*(\Omega^q_X(\log D))\otimes\pi^*\hol(\nu)$
has no global sections, which means that $R^j\pi_*(\Omega^q_X(\log D))$
must be zero.
\qed

Let $\nu =0$. Let $z_i$ be the projective defining equation for $H_i$. 
We fix $H_{i_\infty}$ with $i_\infty \in I$
as the hyperplane at infinity and we will denote it by $H_{\infty}$. 
Let $x_i= z_i/z_{i_\infty}$ and let $w_i =d \log x_i$ 
be the differential form with a logarithmic pole along $H_i$ 
with residue $1$ and a logarithmic pole along $H_\infty$ with 
residue $-1$.
Let $\omega \in H^0(U, \Omega^1_U)$ 
be given as
\begin{equation}
\label{omega}
\omega = \sum_{i\in I \setminus \{i_\infty\}}\alpha_i\omega_i 
\end{equation}
with $\alpha_i \in \C$.  Then $\omega$ has a logarithmic
pole along $H_\infty$ with residue $a_{i_\infty} = 
-\sum_{i\in I\setminus\{i_\infty\}}\alpha_i$.
Let $\tilde\omega = \pi^*\omega$, where $\pi: X \flecha \pn$ is 
the {\it standard} resolution of $H$ and let $D=\pi^*(H)$. As $H^0(U, \Omega^1_U)$
injects into $H^0(X,\Omega^1_X(\log \pi^{-1}(\arrA)))$, we still denote
$\pi^*\omega_i$ again by $\omega_i$.
The form $\omega$ defines a connection $d+\omega$ on the 
$rank~1$ bundle $\hol$ which, as $d\omega =0$, is {\it
integrable,} i.e. it has zero curvature.  Let $U= X\setminus
D$ and $j:U\flecha X$ the inclusion. Let $\Omega^\bullet_U$ 
be the de~Rham complex with
the differential $\nabla = d+\omega$. We have a local constant
system $V$ over $U$ given as $V= ker(\nabla)$. If for every 
$i\in I$ and for every $l\in \calL$, the residues $\alpha_i$ and 
$\sum_{i\in I_l}\alpha_i$ are not positive integers, then 
the cohomology of the local system $V$ is then given by 
\[
H^p(U,V) = 
\H^p(X, Rj_*V)=\H^p(X, \Omega^\bullet_X(\log (D)), \tilde\nabla)
\]
where $\tilde\nabla = \pi^*\nabla$, see [D1, II.6].
Moreover, by Lemma \ref{lemmaVanishing}, see [D1, II.6]  too, 
we have that
\[
\H^p(X, \Omega^\bullet_X(\log (D)), \tilde\nabla) = 
 H^p(H^0(X, \Omega^\bullet_X(\log (D)), \wedge\omega))
\]
where $ H^p(H^0(X, \Omega^\bullet_X(\log (D)), \wedge\omega))$ 
is the $p$-homology of the complex of 
vector spaces $H^0(X, \Omega^\bullet_X(\log (D)), \wedge\omega).$

Let $A^p = H^0(X, \Omega^p_X(\log (D)))$.
The exterior product by $\omega$ induces the complex of vector 
spaces
\begin{equation}
\label{aomotocomplex}
0 \longrightarrow A^0 \stackrel{\omega}{\longrightarrow} A^1
\stackrel{\omega}{\longrightarrow}  A^2
\stackrel{\omega}{\longrightarrow} \ldots 
\stackrel{\omega}{\longrightarrow}  A^n
\longrightarrow 0 .
\end{equation}
This complex appeared for the first time in [A] and will play a
central role in the following sections.

\begin{Theorem}[EV1, 1.5 and 1.7]
\label{theorem vanishing local systems}
Under the hypothesis of Lemma~\ref{lemmaVanishing},
let $\omega = \sum_{i\in I\setminus\{i_0\}}
\alpha_i\omega_i \in H^0(U, \Omega^1_U)$
and $V= ker(d+\omega)$ be a local system.
If $\alpha_i \not\in \Z$ for every $i\in I$ and if for every 
$L\in \calL$ $\sum_{i\in I_L}\alpha_i \not\in \Z$ then, 
for $U$ affine, we have that
\[
\H^p(X, Rj_*V)=\H^p(X, j_!V)=0  \hbox{\ \ \ \ \ \ \ for\ \ } p\neq n .
\]
\end{Theorem}

\begin{Corollary}
\label{Vanishing cohomology from Aomoto}
Under the hypothesis of Theorem \ref{theorem vanishing local systems} 
one has:
\[
\begin{array}{l}
1)~  H^p(A^\bullet, \wedge \omega)= H^p(U, V)
\\ \ \\
2) \hbox{ The complex }(A^\bullet, \wedge \omega) 
\hbox{ is exact, except in degree $n$.}
\end{array}
\]
\end{Corollary}

Let $\Gamma^p \subset A^p$ 
be given as
\[
\Gamma^p = 
\{\bigwedge_{j=1}^p \omega_{i_j} ~|~ i_j \in 
I\setminus \{i_\infty\}\}
\]
where, as above, $\omega_i$ is the pull back of the logarithmic 
differential form $\omega_i$.

We have the following result due to Brieskorn, see [B, Lemma 5]. As
we later, in Section 4, apply a similar method we give a proof here. 

\begin{Claim}
\label{generator of A^p}
The set $\Gamma^p$ generates $A^p$ as a $\C$ vector space. 
\end{Claim}

\prf
As in Lemma \ref{lemmaVanishing} the proof will be by induction
on $|I|$.
For $|I|=1$ we only have one hyperplane, namely the one at infinity
so $\Gamma^p= \emptyset.$
On the other hand, from (\ref{eq:diff.form.seq})
we have that  
\[
H^0(\pn, \Omega^p_{\P^n}(\log (H_\infty))) = 0
\]
for $p>0$.

Let $|I|>1$. 
For $I' \subset I$ a proper not empty subset, we can assume that
$i_\infty \in I^\prime$ otherwise we can choose another hyperplane as
the one at infinity.  Let 
\[
A^{\prime^p}=
H^0(X, \Omega^p_{X}(\log \gamma^*(D')), \wedge\omega^\prime)
\]
where $\pi ' : X' \flecha \pn$ is the standard resolution of 
$H' = \sum_{i\in I'}H_i$,  $D'= \pi'^*(H')$ and $\gamma : X \flecha X'$
is the morphism given by Lemma \ref{lem:blup}.(a) . 
Let $\Gamma^{\prime^p}=\{\bigwedge_{i_j}^p \omega_{i_j} 
~|~ i_j \in I'\setminus \{i_\infty\}\}$. As inductions hypothesis
we assume that the claim holds true 
for any proper subset $I'\subset I$. 
We fix $i_0 \in I$ with $i_0 \neq i_\infty$ and let 
$I'= I\setminus\{i_0\}.$  From (\ref{suc.pullback-formasdif.1})
we have the following exact sequence
\begin{equation}
\label{extension for A^p}
\begin{array}{c}
0\flecha H^0(X, \Omega^p_X(\log \gamma^*(D'))) \flecha 
   H^0(X, \Omega^p_X(\log (D))) 
\\ \ \\
\flecha H^0(D_{i_0}, \Omega_{D_{i_0}}^{p-1}(\log (D')|_{D_{i_0}}))\flecha 0.
\end{array}
\end{equation}

The left map in (\ref{extension for A^p}) is given by the natural
inclusion and the right one is given as
\begin{equation}
\omega_{i_1}\wedge\dots\wedge\omega_{i_q} \mapsto
\left\{\begin{array}{ccc}
0 & \hbox{ \ if \ } &    i_j \neq  i_0 
                        \hbox{\ for \ } 1\leq j \leq p
\\ \ \\
\omega_{i_1}\wedge\dots\wedge
    \widehat{\omega_{i_j}}\wedge\dots\wedge\omega_{i_q}|_{D_{i_0}} & 
               \hbox{\  if \ } & i_j = i_0
                    \hbox{\ for \ } 1\leq j \leq p   
\end{array}\right. .
\end{equation}

By induction on the dimension the restriction of this map to
$A^{\prime^{n-1}}\wedge \omega_{i_0}$ is surjective and one 
obtains
\begin{equation}
\label{descomposicion de A^p}
A^p = A^{\prime^p}+A^{\prime^{p-1}}\wedge \omega_{i_0}.
\end{equation}
Induction on $|I|$ implies the claim.
\qed

\begin{Remark}
\label{descomposicion de A^p como suma directa}
The sum (\ref{descomposicion de A^p}) is a direct sum for the case 
when $l \cap H_{i_0} \neq l $
for every $l \in \calL$, i.e. when $H_{i_0}$ does not contain
``bad loci''. 
\end{Remark}

\prf
The result follows directly from \ref{prop:rel}.
\qed


\section{The Gau\ss-Manin Connection}
\label{section GM-connection}

Let $\arrA$ be an arrangement of $m= n+r+1$ hyperplanes in $\P^n$ which 
does not necessarily have normal crossings. 
We fix once and for all an order
``$<$'' on the set of hyperplanes such that $\{H_0,...,H_n\}$, the first
$n+1$ hyperplanes, are linearly independent.  Let $(z_0,...,z_n)$ be 
homogeneous coordinates for $\pn$. We choose the coordinates in such 
a way that $z_i$
is a homogeneous defining equation of $H_i$ for $i=0,...,n$. We have that
\[
H_j := \sum_{i=0}^n \lambda _{j,i} z_i = 0  \hskip 1cm
(j=n+1,...,m).
\] 

We denote by $ {\P}^{n^\vee}$ the projective space dual to
${\P}^n$.  As every point $p \in {\P}^{n^\vee}$ represents a
hyperplane $H_p \subset {\P}^n$ we can now consider the locus 
in $\pndual$ defined as the set $p \in \pndual$ such that the 
configuration $\arrA \cup H_p \subset \P^{n}$  has more non-normal 
crossings than $\arrA$. This set is known as the discriminant of 
$\arrA$, it forms a divisor in $\pndual$    
and will be denoted by $Discr(\arrA)$. Let 
$S= \pndual \setminus Discr(\arrA)$.  The discriminant is
not necessarily a normal crossings divisor even though if
$\arrA$ was one.

Let $\{h_0,...,h_n\}$ be homogeneous coordinates of $ {\P}^{n^\vee}$
dual to $(z_0,...,z_n)$ and let $\cal S^\vee$ be the homogeneous 
coordinate ring of $ {\P}^{n^\vee}$. We can identify $\cal S^\vee$ 
with the set of homogeneous polynomials in the $h_i$'s. Let 
$J \in M_{m+1,n+1}(\cal S^\vee)$ be given by
\begin{equation}
\label{jacobian}
J = \left( \begin{array}{ccc}
          1          &  \ldots  &      0        \\
          \vdots     &  \ddots  &    \vdots     \\
           0         &  \ldots  &       1       \\ 
       \lambda_{0,0} &  \ldots  & \lambda_{0,n} \\
        \vdots       &   \ddots & \vdots        \\
       \lambda_{r,0} &  \ldots  & \lambda_{r,n} \\
         h_0         &  \ldots  &   h_n
        \end{array} \right)  .
\end{equation}
 
From $J$ one can write the discriminant of ${\arrA}$ as the union of
the zero set of all non trivial $(n+1)$-minors that contain the row
$\{h_0,...,h_n\}$. 
The fact that the arrangement can have non-normal crossings
means that some of the $n+1$ minors of the matrix $J'$ can be
zero. 

Let us consider now a family of arrangements in $\P^n$ given by 
the projection
\[
\pi : \P^{n^\vee} \times \pn \longrightarrow \P^{n^\vee}.
\]
The fiber of $\pi$ on $p=(h_0,...,h_n) \in \pndual$ represents
again the arrangement $\arrA$ plus a hyperplane that moves 
smoothly when we move smoothly on $ \P^{n ^\vee}$.
The extra hyperplane is defined by the fibers of
$\Delta = \{ (h_0,...h_n)\times(z_0,\dots,z_n) \in \P^{n^\vee} 
\times \pn~|~z_h:=\sum_{i=0}^nh_iz_i=0\}$. We have a divisor in 
$\pndual \times \pn$ given as 
$[\pndual \times \arrA] \cup \Delta \cup [Discr(\arrA)\times\pn]$. 
This divisor does not have normal crossings. 
We will denote this divisor by $D$. 

Let 
\[
\varrho : \tilde S \flecha \pndual
\]
be a blow up along the elements of $\calL(Discr(\arrA))$ 
as shown in section \ref{sec. Vanishing Theorems}.  
Let $X = {\tilde S}\times\pn$ and $D^\prime = 
(\varrho\times id_{\P^n})^{-1} D$.
Let $\rho :\tilde X \flecha X$ be a blow up, again as shown in section 
\ref{sec. Vanishing Theorems}, along the elements of $\calL(D^\prime)$
and let $\tilde D = \rho ^* D^\prime$, such that $\tilde D$ has 
normal crossings.  
We have the following diagram
\begin{equation}
\label{family normalization}
\begin{array}{ccccc}
       \tilde X           &         &     &      &\\
\rho \downarrow \hbox{ \  } &\stackrel{~~~~\sigma}{~\searrow}   &&& \\
     X &\flecha & \pndual\times\pn &\stackrel{\pi^\prime}\flecha & \pn\\
   \downarrow     &          & \downarrow \pi\\
        \tilde S       &\stackrel{\varrho}{\flecha}& \pndual.
\end{array}          
\end{equation}
We denote the morphism $\tilde X \stackrel{\rho}{\flecha} X
\flecha \tilde S$ again by $\pi$. 
Let $\Omega^i_{\tilde  X/\tilde S}(\log \tilde D)$ be the coherent sheaf
of $\hol_{\tilde X}$-modules of relative $i$-forms of $\tilde X$ 
over $\tilde S$ with logarithmic poles along $\tilde D$. 

We fix once and for all the hyperplane $H_0$ as the hyperplane at
infinity. For $1\leq i\leq m$ let $\omega_i= d_{rel}x_i/x_i$, with
$x_i = z_i/z_0$,  be the differential
form holomorphic on $U= \pn \setminus \arrA$ with a logarithmic 
pole along $H_i$ with $1$ as residue, and a logarithmic pole along
$H_0$ with residue $-1$.  Let $\tilde\omega_i = 
\sigma^*\pi^{\prime^*}\omega_i$.
Let $\tilde U = \pndual\times\pn \setminus D$ and let $\omega_s = 
d_{rel}x_s/x_s$
be the differential form holomorphic on $\tilde U$ 
where $x_s := 1+\sum_{i=1}^n l_ix_i$ with $l_i=h_i/h_0$,
$x_i=z_i/z_0$ and where the 
differential $d_{rel}$ is the relative differential, i.e. 
$d_{rel}|_{\pi^{-1}\hol_S}\equiv 0$. The differential form $\omega_s$ has then
a logarithmic pole along $\Delta$ with residue $1$, and a logarithmic 
pole with residue $-1$ along $z_0 = 0$.

\begin{Remark}
\label{residue for omega_s with absolute differential}
The absolute differential form
$\omega_s= dx_s/x_s$ has also a logarithmic pole along $h_0=0$ with
residue $-1$.
\end{Remark} 

\begin{Notation} 
\label{notation ueber die Differentiale Formen} 
Let $H_0$ and $H^\prime$ two
hyperplanes in $\pn$ with $z_0$ and $z^\prime$ as homogeneous defining
equations. We take $H_0$ as the hyperplane at
infinity. Let $x^\prime = z^\prime/z_0$ be the affine defining
equation of the hyperplane induced by $H^\prime$ on the affine space
$\pn \setminus H_0$. We denote by $[\frac{dx^\prime}{x^\prime}]$ or
by $[\frac{dz^\prime}{z^\prime}-\frac{dz_0}{z_0}]$ the global 
differential form, holomorphic on $U = \pn \setminus H^\prime \cap H_0$, with
a logarithmic pole along $H^\prime$ with residue $1$ and a logarithmic
pole along $H_0$ with residue $-1$.
\end{Notation}

Let $\omega \in H^0(\tilde U, \Omega_{\tilde U/S}^1)$
be given by
\[
\omega = \sum_{i=1}^m a_i\omega_i + a_h\omega_s.
\]  
The pull back $ \sigma^*\omega$ induces the differential form  
$\tilde\omega \in  H^0(\tilde X, 
   \Omega^1_{\tilde X/\tilde S}(\log \tilde D))$.  This form
has residue $a_i$ along the pole $H_i$, $a_h$ as residue along
$\Delta$ and $\sum_{i\in I_L} a_i$ as residue along the exceptional
divisor $e_L = \sigma^{-1} (L)$ for $L \in \calL(D)$ .

We consider the operator $\nabla = d_{rel}+ \tilde\omega$. 
As $d_{rel}\tilde\omega = 0$ we
have a logarithmic de~Rham complex 
\begin{equation}
\label{complxrel GM}
0\flecha\hol_{\tilde X}\flecha\Omega^1_{\tilde X/\tilde S}(\log\tilde D)
\flecha\dots\flecha\Omega^n_{\tilde X/\tilde S}(\log\tilde D)\flecha 0.
\end{equation}

The $ker(\nabla:\hol_{\tilde X}\flecha\Omega^1_{\tilde X/\tilde S}
(\log\tilde D))$ defines 
a relative local constant system $V_{rel}$ over the complement of $\tilde D$
in $\tilde X$ relative to $\tilde S$. 

Let $H^i_{DR}( \tilde X/\tilde S, \tilde D, \nabla)$ be the $i$-th 
de~Rham cohomology of $ \tilde X$ relative to $\tilde S$ with 
respect to $\nabla$. This is the sheaf of $\hol_{\tilde S}$ modules 
\[
H^i_{DR}(\tilde X/\tilde S, \tilde D, \nabla) = {\bf R}^i\pi_*
(\Omega^\bullet_{\tilde X/\tilde S}(\log {\tilde D}), \nabla)
\]
where ${\bf R}^i\pi_*\Omega^\bullet_{\tilde X/\tilde S}(\log {\tilde D})$ 
are the hyperderived functors of the functor ${\bf R}^0\pi_*$.
Under the assumptions that all the residues of $\omega$ are not
positive integers, we have, see [D1, II.6], that these are the sheaves 
of cohomology groups of the relative local system $V_{rel}$.

As $\pi: \tilde X \flecha \tilde S$ we have the following 
exact sequence
\begin{equation}
\label{differential extension}
0\flecha \pi^*\Omega^1_{\tilde S}(\log (\varrho^*Discr(\arrA))) 
    \flecha\Omega^1_{\tilde X}(\log\tilde D) \flecha
        \Omega^1_{\tilde X/\tilde S}(\log\tilde D)\flecha 0.
\end{equation}

One can extend the differential on the relative complex by taking
$\tilde \Omega \in \newline
H^0(\tilde X, \Omega^1_{\tilde X}(\log\tilde D))$ 
as $\sigma^*(\sum^m_{i=1} a_i\frac{dx_i}{x_i} + a_h\frac{dx_s}{x_s})$
where again $x_i$ and $x_s$ are as before and where $d$ is the absolute 
differential.
We define $\bar \nabla = d + \tilde\Omega$. As 
$d\tilde\Omega=0$ we have $\bar\nabla^2=0$. This defines
the complex $(\Omega_{\tilde X}^\bullet(\log\tilde D), \bar\nabla)$.
Let $V=ker(\nabla:\hol_{\tilde X}\flecha\Omega^1_{\tilde X/\tilde S}
(\log\tilde D))$.  Filtering the complex 
$(\Omega_{\tilde X}^\bullet(\log\tilde D), \bar\nabla)$ by 
\[
\dots\bfF^{i+1}\subset\bfF^i\subset\dots
\subset\bfF^0=\Omega_{\tilde X}^\bullet(\log\tilde D)
\]
where
\[
\bfF^i = \pi^*\Omega_{\tilde S}^i(\log (\varrho^*Discr(\arrA)))\wedge
\Omega^{\bullet -i}_{\tilde X}(\log\tilde D)
\]
we can construct a spectral sequence abutting to 
${\bf R}^\bullet\pi_*(\Omega^\cdot_{\tilde X}(\log\tilde D))$. The
$E_1^{a,b}$ terms are equal to 
$\Omega^a_{\tilde S}(\log \varrho^*(Discr(\arrA))) 
\otimes_{\smholS}{\bf R}^b
\pi_*(\Omega^\bullet_{\tilde X/\tilde S}(\log\tilde D))$
and the differential 
\begin{equation}
\label{gauss-manin differential GM}
d_1:E_1^{a,b}\longrightarrow E_1^{a+1,b}
\end{equation}
has
bidegree $(1,0)$.

From the filtration we have
\[
0\flecha \frac{\bfF^1}{\bfF^2}\flecha \frac{\bfF^0}{\bfF^2}
\flecha \frac{\bfF^0}{\bfF^1}\flecha 0.
\]
This is just the exact sequence of complexes
\begin{equation}
\label{exact sequence GM}
\begin{array}{c}
0\flecha \pi^*\Omega^1_{\tilde S}(\log \varrho^*(Discr(\arrA))\otimes 
\Omega^{\bullet -1}_{\tilde X/\tilde S}(\log\tilde D)\\
\\
  \flecha \frac{\Omega_{\tilde X}^\bullet(\log\tilde D)}
     {\pi^*\Omega_{\tilde S}^2(\log \varrho^*Discr{\small(\arrA)})
      \wedge\Omega^{\bullet -2}_{\tilde X}(\log\tilde D)}
\flecha \Omega^{\bullet}_{\tilde X/\tilde S}(\log\tilde D) \flecha 0.  
\end{array}
\end{equation}
The differential (\ref{gauss-manin differential GM}), for the case
when $a=0$, is the connecting morphism for the long exact sequence
of cohomology obtained from (\ref{exact sequence GM}). Using
projection formula, with local calculation one can show that it has
the Leibniz properties of a connection.   
It is called the Gau\ss-Manin connection, see [K,4.6]. We will 
still denote it by ${\bar\nabla}$.

For the integrability of the Gau\ss-Manin connection we have the
following diagram.
\begin{equation}
\label{integrabilitaetsdiagram}
\begin{array}{ccccccccc}
  &      &   0         &        &    0      &        &   0        &  &\\
  &      &\uparrow     &        &\uparrow   &        &\uparrow    &  &\\
0 &\flecha &\bf{F^1/F^2} &\flecha &\bf{F^0/F^2} &\flecha &\bf{F^0/F^1} 
                                                           &\flecha &0\\
&&&&&&&&\\
  &        &\uparrow   &        &\uparrow  &         &     ||     & &\\
&&&&&&&&\\
0 &\flecha &\bf{F^1/F^3} &\flecha &\bf{F^0/F^3} &\flecha &\bf{F^0/F^1} 
                                                            &\flecha &0\\
&&&&&&&&\\
  &        &\uparrow &        &\uparrow &        &      &        &\\
&&&&&&&&\\
  &        &\bf{F^2/F^3} &=     &\bf{F^2/F^3} &    &        &        &\\

  &        &\uparrow &        &\uparrow &        &          &     &\\
  &        &    0    &        &    0    &        &          &     &.
\end{array}
\end{equation}
The curvature is then given by the map
\[
\bar\nabla^2: {\bf R}^a\pi_*(\bf{F^0/F^1})\flecha 
{\bf R}^{a+2}\pi_*(\bf{F^2/F^3}).
\]
For an element $\alpha \in {\bf R}^a\pi_*(\bf{F^0/F^1})$, the connecting
morphism of the middle horizontal exact sequence in 
(\ref{integrabilitaetsdiagram}) gives us an element in 
${\bf R}^{a+1}\pi_*(\bf{F^1/F^3})$.  From the 
left vertical exact sequence
in (\ref{integrabilitaetsdiagram}) one has the integrability of the
Gau\ss-Manin connection.

\begin{Proposition}
\label{prop. existence of GM-connection}
For $\pi$ and $\omega$ as above we have that
\\ \
\hbox{\hskip12pt}
(i) $H^i_{DR}(\tilde X/\tilde S, \tilde D, \nabla)$ has an
integrable connection.
\\ \
\hbox{\hskip12pt}
(ii) Under the assumptions that the residues of $\omega$ are
not integers then 
${\bf R}^i\pi_*(U, V)$ is a local system 
equal to
${\bf R}^i\pi_*
(\Omega^\bullet_{\tilde X/\tilde S}(\log {\tilde D}), \nabla)$
where $V$ is the local system of flat sections of $\bar\nabla$.

\end{Proposition}


\section{Some Combinatorics}
\label{section combinatorics}

Let $\arrA$ be an arrangement of $m=n+r$ hyperplanes in the affine space 
$\C^n$. 

\begin{Definition}
Let $L=L(\arrA)$ be the set of all non-empty   
intersections $\cap_{i\in I}H_i \not=\emptyset$ of elements of $\arrA$.  
We assume $\C^n \in L(\arrA)$ as the intersection over the empty set, 
i.e. $I=\emptyset$. 
\end{Definition}

One has a partial 
order ``$\preceq $'' on  the elements of $L$  given by reverse inclusion:
$X \preceq Y$, for $X, Y \in L$ if and only if $Y
\subseteq X$.  
Let $X,~Y \in L$ be such that $X\prec Y$. A $chain$ from $X$ to $Y$
is a set $\{ Z_0,\dots,Z_n\}\subset L$
such that $X=Z_0\prec Z_1\prec\dots\prec Z_n=Y$.  We say then, that 
this chain has $length$ $n$.  A chain is $maximal$ if for every
$i \in \{0,\dots,n-1\}$ there does not exist $W\in L$ such that
$Z_i\prec W\prec Z_{i+1}$. 

We have two binary operations on the elements of $L$, namely the $meet$
defined as $X \wedge Y = \cap \{Z \in L | X \cup Y \subseteq Z\}$ and the
$join$ defined as $X \vee Y = X \cap Y$ when $X\cap Y \not=\emptyset$.  
We have a $rank$ function given
as $rank(X) = codim X$ for $X \in L$. For the case when $X\vee Y$ exists,
this function satisfies 
\begin{equation}
\label{rank inequality}
rank(X\wedge Y)+rank(X\vee Y) \le rank(X) +rank(Y).
\end{equation}

\begin{Lemma}
\label{lemma length maximal chains}
If $X,~Y \in L$ and $X\prec Y$, then any maximal  
chain from $X$ to $Y$ has the same length and one says that
$L$ satisfies the Jordan-Dedekind chain condition.
\end{Lemma}

\prf The lemma follows from seeing that the length of any 
maximal chain is equal to $rank(Y)-rank(X)$. 
\qed

\begin{Definition}
Given $X, Y \in L$ such that $X \preceq Y$ we define the closed
interval as $[X,Y] = \{Z\in L ~|~X \preceq Z\preceq Y\}$. 
One can define the open interval by
taking strict inequalities.
For $X\in L$ we say it is {\bf maximal} (resp. {\bf minimal})
if there does not exist $Y\in L$ such that $X\prec Y$ $(resp.~Y\prec X)$.
\end{Definition}

A partially ordered set $L$ has the structure of a lattice if
the two operations ``$\wedge$'' and ``$\join$'' exist
for every pair of elements in $L$. A lattice is called geometric
if there is a rank function $rank: L \flecha \N$ satisfying 
(\ref{rank inequality}).
If $L$ is a partially ordered set for which only $X\wedge Y$ exists for all
$X,Y \in L$ then $L$ is called semi-lattice.

One can see that $L(\arrA)$ has the structure of a geometric 
semi-lattice and is called the $intersection~semi$-$lattice$ of 
the arrangement $\arrA$. 

If $\arrA$ is central, (i.e. $\cap_{H \in \arrA} H    
\neq \emptyset$) then the ``join'' of any two elements exists and $L(\arrA)$ 
has a unique maximal element. We have then that $L(\arrA)$ has the
structure of a geometric lattice.  

From a $n$-dimensional non-central arrangement one can construct a 
$(n+1)$-dimensional central arrangement. This process is know as the
coning process and is given as follows:\newline
Let $f_i \in \C[x_1,\dots,x_n]$ the defining equation for the
hyperplane $H_i \in \arrA$ with $\arrA$ non-central.   Let 
$F_i \in \C[z_0,\dots,z_0]$ be the homogenized polynomial of $f_i$ obtained
by substituting $x_j=z_j/z_0$ in  $f_i$ and multiplying by $z_0$.
Let $F_0= z_0$. We define the $cone$ over $\arrA$ as
the $(n+1)$-dimensional arrangement $\coneA= \cup_{i=0}^n \bar H_i$  where
$\bar H_i = ker(F_i)$.

As the arrangement $\coneA$ is central, the inequality 
(\ref{rank inequality}) is always satisfied.  
Elements $l\in L(\arrA)$ of 
codimension $n-r$ (resp. dimension $r$) give rise to elements of 
codimension $(n+1)-r$ (resp. dimension $r$) of $L(\coneA)$. 
We have that $|\coneA| = |\arrA| + 1$.

The deconing process is then given as the inverse process of the coning
process.
Let $\arrA = \cup_{i\in I}\bar H_i$ be a central arrangement. We fix 
a hyperplane $\bar H_0$ as the hyperplane
at infinity.  Let $F_i\in \C[z_0,\dots,z_n]$ be the homogeneous 
defining equations for $\bar H_i$. We can choose the set of coordinates
such that $H_0$ is given by $z_0=0$. Let then $\arrA^\prime$ be the affine
arrangement given by the hyperplanes $H_i$ with defining equation
$f_i \in \C[x_1,...,x_n]$ such that $f_i=F_i/z_0$ and where $x_i = z_i/z_0$.

With this method one can translate the definitions and results from
central arrangement to non-central arrangements and back.

\begin{Definition}
\label{definition levels}
Let $L_{n-r}=L^r= \{Z \in L~|~ rank(Z) = r \}$. 
We have that $L^0 = \C^n$, $L_{n-1}=L^1$ is the set of hyperplanes 
of our arrangement and $L_0$ is a set of points in $\{\C^n\}$.  
\end{Definition}

\begin{Lemma} 
Maximal elements of $L(\arrA)$ have the same rank.
\end{Lemma}

\prf See [OT, lemma 2.4]
\qed

We have the following definition.

\begin{Definition}
\label{definition rank of a lattice}
We define $rank$ of $L=L(\arrA)$ as the rank of any maximal element.
\end{Definition}  

Let ``$<$'' be an arbitrary but fixed linear order for the elements 
of $\arrA$, i.e. we fix a linear order for the elements of $L^1$.
Let ${\cal E}_1$ be the complex vector space freely generated by elements
$\{ e_H : H \in \arrA \}$ and let $\cal E$ be its exterior algebra. For 
$S \subset \arrA $ we denote the element $\wedge_{H\in S }e_H$
by $e_S$ respecting the order chosen. A subset $S \subset \arrA$ 
is said to be $dependent$ if there exists $H_0 \in \arrA$ such that
$\cap_{H\in S \backslash \{H_0\}}H = \cap_{H\in S}H$.
We have the morphism 
\[
\partial : {\cal E} \longrightarrow {\cal E}
\]
where for $e_{(H_1,...,H_p)} \in {\cal E}_p$ is given as
$\partial e_{(H_1,...,H_p)} = \sum_{i=1}^p (-1)^i e_1\wedge...\wedge 
\widehat{e_i}\wedge...\wedge e_p$.  It is a morphism of algebras
and  it is easy to see that $\partial^2=0$.
Let  $\cal J$ be the graded ideal of $\cal E$ generated by 
$\partial e_S$ with $S\in \arrA$ dependent. 
The quotient $\cal E/J$ is a graded algebra known as the Orlik-Solomon 
algebra and appeared
for the first time in [OS], see [OT] too. 
From [OS, Theorem 5.2] we have the isomorphism 
\begin{equation}
\label{Orlik-Solomon presentation}
{\cal E/J} \flecha H^*(U, \C)
\end{equation}
with $U=\C^n \backslash \arrA $. 

For the intersection semi-lattice $L(\arrA)$ we have the following 
definition of a neat base-family which we will denote as Nbf. 
As before we denote by $0$ the minimal element of our semi-lattice $L$. 
\begin{Definition}
\label{neat base family}
(i) If $rk(L) = 1$, then $\{H \in L ~|~ H~maximal \}$ is a Nbf.

(ii) We assume the existence of Nbf for lattices of
rank$\le (n-1)$ and suppose that $rk(L)=n$. For every $p \in L^n$, an
upper bound of $L$, 
chose a fixed $H_p \in L^1$.
For every $l \in L^{n-1}$ such that $H \not\preceq l\prec p$ let 
${\bf B}_l$ be a Nbf for the lattice $[0,l]$.  We define 
\[
{\bf B} = \cup_{p\in L^n}\{B\cup \{H_p\} / 
B\in \cup_{H\not\preceq l} {\bf B}_l\}
\]
is a Nbf for $L$.
\end{Definition} 

This definition was first made by A.~Bj\"orner for the case when
$L$ is a lattice. With respect to his definition the neat base 
family defined above is the union of the neat base families of
the lattices $[0,l]$ with $l\in L(\arrA)$ maximal, see [Bj, Section 3].

Given a subset $S \subset L^1$  of $points$ in a semi-lattice we have that
$rank(\vee_{H\in S} H) \leq \newline 
\mid S \mid$.  We have the following 
definition, see [OT, definition 3.2].

\begin{Definition}
\label{definition dependent set}
For $S\subset L^1$, with $L^1$ as in definition \ref{definition levels}
we say that $S$ is {\bf independent} if when $\join_{H\in S}H$ exists
$rank(\vee_{H\in S} H) = \mid S \mid$
otherwise we say that $S$ is {\bf dependent}.
\end{Definition}

This notion of dependence is compatible with the one above.  

\begin{Definition}
Maximal independent subsets will be called {\bf bases} and the minimal
dependent sets will be called {\bf circuits}. 
\end{Definition}

We will denote the set of
independent elements of $L$ by ${\cal I}(\arrA)$ and the set of 
circuits by ${\cal C}(\arrA)$.  In our case, one can easily see 
that, for ${\cal L}$ as in Definition \ref{def. bad}, we have that 
$l \in {\cal L}$ if and only if the set $\{ H \in I_l \}$ contains a 
circuit. 

\begin{Definition} 
Let $C \subset L^1$ be a circuit and
$p\in L^1$ the least element of $C$.  Then we say that $C - \{p\}$ 
is a {\bf broken circuit}. For a broken circuit $C \subset L^1$ 
let $princ(C) \in L^1$ be such that $C \cup \{princ(C)\}$ is a circuit
and is the smallest element with this property.  The family of sets 
which do not  contain
a broken circuit will be called $non\ broken\ circuits$; we will call
them nbc-elements. A maximal nbc-element will be called an nbc-base. 
\end{Definition}

We have the next proposition due to A. Bj\"orner, see [Bj, Proposition 3.8]

\begin{Proposition}
\label{proposition nbc are Nbf}
Let $L$ be a geometric semi-lattice of rank r and ``$\preceq$'' a linear
order of the elements of $L^1$. Then the collection nbc-bases are
an Nbf for $L$. 
\end{Proposition}

Let $x_i$ be the defining equation for the hyperplane $H_i \in \arrA$
and let  
\[
A^p=\{\wedge_{j=1}^p \frac{dx_{i_j}}{x_{i_j}}~|~ i_j \in {1,\dots,m}\}. 
\]
One can obtain 
a basis for $A^n$ and $A^{n-1}$ by following the
construction defined in [SV, 1.6].  To each element
$X\in L$ we associate a hyperplane $H_X$ such that $H_X \preceq
X$ and such that for every $H < H_X$ we have that $H \not\preceq X$. 
Taking complete flags over $L$ of the form 
$X_{l-1}\prec...\prec X_1\prec X_0$
where $X_i \in L_i$ and is such that $H_X\not\prec X_i$ when $X$ is an 
element of the flag and $X_i\prec X$.  The set of such flags of 
length $l$ will be denoted by $FL_l$. To every such flag 
$X_{n-1}\prec...\prec X_1\prec X_0 \in Fl_n$ one can associate
the $n$-form $\omega_{\alpha_0}\wedge...\wedge \omega_{\alpha_{n-1}}
\in A^n$, where $\alpha_i$ is a defining equation for $H_{X_i}$ and
$\omega_{\alpha_i}=d{\alpha_i}/{\alpha_i}$. From construction is 
clear that our $n$-form is not trivial.  We denote by
$B^n$ the set of all $n$-forms obtained in this way. Following
the same construction but now taking flags starting over elements $X
\in L_{n-r}$. In this case, let us denote the complete flags of
length $l\leq r$ by $FL_l^r$ and by $B^r$ the set of all
the $r$-forms obtained from elements in $FL_r^r$, that is, the
forms that come from maximal complete flag constructed over
elements of our intersection lattice of codimension $r$. 
One can see that $FL_r^r$ forms a neat base family for the lattice 
$L(\arrA)$.  We have the following theorem due to A. Bj\"orner, see 
[Bj, Theorem 4.2] and [SV].

\begin{Theorem}[Bj] 
\label{teo:nbc}
The set of $B^r$ forms a basis for $A^r$.
\end{Theorem}

We would like to calculate a basis for $ker (\phi:A^n\flecha H^n(U,\C))$.
By theorem \ref{teo:nbc} we know that the set of $nbc$-bases
is a basis for $A^n$.
Let $\arrA=\cup_{i=1}^{m=n+r}H_i \subset \C^n$ be an arrangement,
${\cal E}_1$ be the complex vector space generated by elements
$\{ e_H | H \in \arrA \}$ and $\cal E/J$ be the Orlik-Solomon algebra.
Let $B \subset \arrA$ a base but not an $nbc$, then $B$ contains at least
one broken circuit. Take $\hat H \in \arrA$ to be the least
element in $\arrA$ with the property that there exists $ C \subset B$ such
that $C$ is a broken circuit and $\hat H = princ(C)$, i.e. $\hat H$ is the 
least element in $\arrA$ such that $\{H\}\cup C$ is dependent.  Let 
$\hat C$ denote the circuit $C\cup \hat H$ and $\hat{B} = B \cup {\hat H}$.

Let ${\cal J}_n = {\cal J} \cap {\cal E}_n$ where $\cal J$ and
${\cal E}_n$ are as before.  From the isomorphism 
(\ref{Orlik-Solomon presentation}), to give a basis 
for $ker(\phi: A^n \flecha H^n(U,\C))$ is equivalent to give a 
basis for ${\cal J}_n$.  We have the following
proposition obtained together with V.~Welker.

\begin{Proposition}
\label{prop:rel}
Let $\arrA=\cup_{i=0}^{m=n+r}H_i \subset \C^n$ be an arrangement,
$\cal E$ the Orlik-Solomon algebra, then a basis for ${\cal J}_n$ is 
given by elements of ${\cal E}_n$ of the form:

i) $e_B$ for $B \subset \arrA$ dependent and $\mid B \mid = n$.

ii) $\partial e_{\hat{B}}$ for $B$ a base but not an nbc.
\end{Proposition}

\prf
We have that
\[
|A^n| = |nbc\hbox{-bases}| +
     \{|\hbox{circuits}|+|\hbox{broken circuits}|\}.
\]
As the number of elements in $(i)$ together with $(ii)$ equal 
$|A^n| - |nbc|$, to show that they form a basis we only have
to show that they generate $ker(\phi)$.
We will prove this by induction by showing that, with help of the
elements in $(i)$ and $(ii)$ one can write any element of $A^n$
as a linear combination of non broken circuits.
Taking the lexicographic order on the elements of $A^n$, induced
by the order chosen for the set of hyperplanes, the induction 
will be taken on this order.

As the first element of $A^n$ is already a non broken circuit so
the statement is true for the base of induction.  Let $C=(H_1,\dots,H_n)$ 
be a base but not nbc.  By inductions hypothesis we can assume that
the statement is true for any base $B\leq C$.  Let $\hat C = C \cup
princ(C)$ as before, then
\begin{equation}
\label{lexicographic induction}
\partial e_{\hat C} = \sum_{i=1}^{princ(C)}e_{\hat C -\{H_i\}}
                   + \sum_{i=princ(C)}^n e_{\hat C -\{H_i\}}.
\end{equation}
Every summand in the first sum of (\ref{lexicographic induction}) contains
a circuit so they all are elements of $i$.  For the second sum, the first
element is $C$ and the rest contain $princ(C)$ so they have 
lexicographic order smaller than $C$.
Applying our induction hypothesis we can write $C$ as linear combination
of circuits and $nbc$-bases.
\qed

\begin{Remark}
One can obtain a basis for the graded ideal $\calJ$ of relations
for the Orlik-Solomon Algebra by defining
a basis for every level $\calJ_r$ in a similar way as in proposition
\ref{prop:rel}.
\end{Remark}

\noindent{\bf Example}

Let $\arrA = \cup_{i=0}^5H_i$ be an arrangement in $\P^2$
given by:
\begin{equation}
\label{projective Ceva}
\begin{array}{l}
H_0 := z_0 = 0 \\
H_1 := z_1 = 0 \\
H_2 := z_1 = 0 \\
H_3 := z_3:= z_0-z_1=0\\
H_4 := z_4:= z_0-z_2=0\\
H_5 := z_5:= z_1-z_2=0
\end{array}
\end{equation}
\[
\begin{array}{c}
\unitlength1cm
\begin{picture}(7,5)
\put (1.7,0.7){\line(1,1){2.7}}
\put (2.3,1.9){$\small{H_5}$}
\put (1.4,3.3){\line(2,-1){5}}
\put (2.3,2.9){$\small{H_4}$}
\put (1.8,5.4){\line(1,-2){2.4}}
\put (3.2,1.5) {$\small{H_3}$}
\put (1.5,1){\line(1,0){5}}        
\put (5, 0.5) {$H_2$}
\put (2,0.4){\line(0,1){5}}        
\put (1, 2 ) {$\tiny{H_1}$}
\put (1.7,5.3){\line(1,-1){4.6}}
\put (4.7,2.5) {$H_0$}
\put (3, 0){$Fig.~1$} 
\end{picture}
\end{array}
\]

Take $H_0$ to be the hyperplane at infinity.  On the
affine complement of $H_0$ 
we have the following arrangement.

\begin{equation}
\label{affine Ceva}
\begin{array}{l}
L_1:= x_1= 0\\
L_2:= x_2= 0\\
L_3:= x_3:= x_1-1=0\\
L_4:= x_4:= x_2-1=0\\
L_5:= x_5:= x_1-x_2=0
\end{array}
\end{equation}
\[
\begin{array}{c}
\unitlength1cm
\begin{picture}(6,5)
\put (1.7,1.7){\line(1,1){2.9}}
\put (2.6,3.2){$L_5$}
\put (1.5,4.3){\line(1,0){3.2}}
\put (2.7,4.5){$L_4$}
\put (4.3,1.5){\line(0,1){3.2}}
\put (4.5,3){$L_3$}
\put (1.5,2){\line(1,0){3.2}}        
\put (2.7,1.6){$L_2$}
\put (2,1.5){\line(0,1){3.2}}        
\put (1.4,3){$L_1$}
\put (2.3,0.6){$Fig.~2$}
\end{picture}
\end{array}
\]

We have that the set of circuits is
\[
\begin{array}{c}
\calC(\arrA) =\{ (L_1,L_3),(L_2,L_4),\\
\\
 (L_1,L_2,L_5) , (L_3,L_4,L_5) \}
\end{array}
\]
and are the only $dependent$ subsets of $\arrA$.
The $nbc$'s are
\[
nbc(\arrA)= \{ (L_1,L_2), (L_1,L_4), (L_1,L_5) ,
	    (L_2,L_3), (L_3,L_4), (L_3,L_5) \}.
\]
Clearly we have that the $broken~circuits$ are only 
$\{ (L_2,L_5), (L_4,L_5) \}$ for which $(L_1) = princ(L_2,L_5)$
and $(L_3) = princ(L_4,L_5)$.

Let $\calE_1$ be freely generated by $\{ e_i | L_i \in \arrA \}$
and let $\calE$ be its exterior algebra.  Let $\calJ$ the ideal
of $\calE$ generated by $\partial e_S$ for $S\subset \arrA$ 
dependent. 
From proposition \ref{prop:rel} we have that $\calJ_n = \calJ \cap
\calE_n$ is generated by
\[
\begin{array}{c}
\calJ_n = <~e_{13}~;~e_{24}~;~ \partial e_{125}~;~ \partial e_{345}~>\\
\\
= <~e_{13}~;~e_{24}~;~e_{12}-e_{15}+e_{25}~;~e_{34}-e_{35}+e_{45}~>.
\end{array}
\]
Under the natural identification of $\calE_1$ with 
$\Omega_U^1$, where $U = \P^2 \backslash \arrA$, given
by $e_i \mapsto \frac{dx_i}{x_i}$, the these relations gives place 
to the following relations of $2$-forms
\begin{equation}
\label{relations ceva projective}
\begin{array}{c}
\frac{dx_1dx_3}{x_1~x_3} =0
\\ \\
\frac{dx_2dx_4}{x_2~x_4} =0
\\ \\
\frac{dx_1dx_2}{x_1x_2}-\frac{dx_1dx_5}{x_1x_5}+\frac{dx_2dx_5}{x_2x_5}=0
\\ \\
\frac{dx_3dx_4}{x_3x_4}-\frac{dx_3dx_5}{x_3x_5}+\frac{dx_4dx_5}{x_4x_5}=0.  
\end{array}
\end{equation}
By proposition \ref{prop:rel} these relations are linearly independent.


\section{The Gau\ss-Manin Matrix}
\label{sec. non normal crossings}

Let $\arrA$ be an arrangement of $m= n+r+1$ hyperplanes in $\P^n$ 
as taken in Section \ref{section GM-connection}. We consider
the same situation as in Section \ref{section GM-connection} but
to make things easier we don't compactify the space of parameters.
We have a family of arrangements in $\P^n$ given by 
the projection
\[
\pi : \XX \longrightarrow S
\]
where $S=\pndual \setminus Discr(\arrA)$.
Let $D$, $U = \XX \setminus D$ and $\omega \in H^0(U, \Omega^1_U)$ be 
given as in Section \ref{section GM-connection} but restricted to $\XX$.

Let $\rho : \tilde X \flecha \XX$ be the blow up, as in 
(\ref{family normalization}) along $\calL (D)=
\calL(S\times\arrA+\Delta\cap[S\times\pn])$ 
taken in the same way as in Section \ref{sec. Vanishing
Theorems}. As our space of parameters is taken as the non compactified
space, we have that 
$\calL(S\times\arrA+\Delta\cap[S\times\pn])= \calL(S\times\arrA)$.
We have the following remark.

\begin{Remark}
\label{remark familie's bad loci}
We have that $\calL(D) = S\times \calL(\arrA) \subset D$, i.e., the
bad loci can only have at most codimension $n$. 
Letting $\pi^\prime : \tilde Y \flecha \pn$ be the
standard resolution along elements of $\calL(\arrA)$
as described in Section \ref{sec. Vanishing Theorems},
one has $\tilde X = S\times \tilde Y$.
\end{Remark}

We denote $\rho^*(D)$ by $\tilde D$.
Let $\Omega^i_{\tilde X/S}(\log \tilde D)$ be the coherent sheaf
of $\hol_{\tilde X}$-modules of relative $i$-forms of $\tilde X$ 
relative to $S$ with logarithmic poles along $\tilde D$. 
Let $\tilde \omega = \rho^*\omega \in 
H^0(\XX, \Omega^1_{\tilde X/S}(\log \tilde D))$. Then $\tilde \omega$
is the differential form with residues $a_h$ along
$\sum_{i=0}^nh_iz_i = 0$, $a_i$ along $H_i$ with 
$a_0 = \sum_{i=1}^m -a_i -a_h$ and where for every $L\in \cal(\arrA)$
the form $\tilde \omega$ has residue $\sum_{i\in I_L}a_i $ 
along the exceptional divisor $e_L = \rho^{-1}(S\times L)$.

We consider lifting of $ \nabla$ to $\tilde X$ given as the operator 
\[
\tilde \nabla = d_{rel}+ \tilde \omega.
\]
Again, as $d_{rel}\tilde \omega = 0$, it gives 
a logarithmic de~Rham complex 
\begin{equation}
\label{complxrel NNC}
0\flecha\hol_{\tilde X}\stackrel{\tilde \nabla}{\flecha}
\Omega^1_{\tilde X/S}(\log \tilde D)
\stackrel{\tilde\nabla}{\flecha}\dots\stackrel{\tilde\nabla}{\flecha}
\Omega^n_{\tilde X/S}(\log \tilde D)\flecha 0.
\end{equation}

The $ker(\tilde \nabla:\rho^*\hol_{\tilde X}\flecha\Omega^1_{\tilde X/S}
(\log \tilde D))$ defines 
a relative local system $V_{rel}$ over the complement of
$\tilde D$
in $\tilde X$ relative to $S$. 

From Proposition \ref{prop. existence of GM-connection} we have
that $H^i_{DR}(\tilde X/S,  \tilde D, \tilde \nabla)$
are the cohomology sheaves of groups of the relative local 
system $V_{rel}$ when the residues of $\tilde\omega_s$ along the components
of $\rho^*(D)$ are not positive integers, see [D1, II.6].

Let $A^p\subset \pi_*\Omega^p_{U/S}$ be 
generated over $\hol_S$ by
\[
\{\wedge_{j=1}^p\frac{dx_{i_j}}{x_{i_j}};  
\wedge_{j=1}^{p-1}\frac{dx_{i_j}}{x_{i_j}}\wedge\frac{d_{rel}x_s}{x_s}
\mid
i_j \in \{1,...,m\}\}.
\]

We have the subcomplex $A^\bullet \subset 
\pi_*\Omega^\bullet_{U/S}$ given by
\begin{equation}
\label{global section complex NNC}
0 \longrightarrow \hol_S
\stackrel{\nabla}{\longrightarrow} A^1
\stackrel{\nabla}{\longrightarrow} \ldots 
\stackrel{\nabla}{\longrightarrow} A^n
\longrightarrow 0.
\end{equation}

For $s\in S$, we consider the restriction $ \nabla_s$ of the 
connection to the fiber $\pi^{-1}(s) = s\times \pn \simeq \pn$. 
From Remark \ref{remark familie's bad loci} the restriction 
$\rho|_s$ of the blow up $\rho$ to the 
fiber $\pi^{-1}(s)$ is a resolution of the configuration 
$\arrA \cup H_s \subset \pn$.  
We will denote $\rho|_s^*(\pi^{-1}(s))$ and 
$\rho|_s^*(\arrA \cup H_s)$ as $\tilde X|_s$ and 
$\tilde D_s$ respectively.  
Let $A^p_s$  be 
the finite dimensional subspace of 
$H^0(\pn\setminus \{\arrA\cup H_s\}, 
\Omega^i_{\pn\setminus \{\arrA\cup H_s\}})$  generated by 
$\{\bigwedge_{j=1}^p \frac{dx_{i_j}}{x_{i_j}} \hbox{;} \
\bigwedge_{j=1}^{p-1}
\frac{dx_{i_j}}{x_{i_j}} \wedge
\frac{d_{rel}x_s}{x_s} / \hbox{\ }i_j \in
\{1,...,m\} \hbox{\ and \ } x_{i_j} < x_{i_k} \hbox{\ if\ } j<k
\}$. We have the following complex:
\begin{equation}
\label{eq:complex|_s NNC}
0 \longrightarrow \C
\stackrel{\nabla_s}{\longrightarrow}  A^1_s
\stackrel{\nabla_s}{\longrightarrow}  A^2_s
\stackrel{\nabla_s}{\longrightarrow} \ldots 
\stackrel{\nabla_s}{\longrightarrow}  A^n_s
\longrightarrow 0 .
\end{equation}

By Lemma~\ref{lemmaVanishing} we have that 
$H^p(\tilde X|_s, \Omega^i_{\tilde X|_s}(\log \tilde D|_s)) = 0$ 
for $p>0$. 
This implies that the cohomology of
the de~Rham complex $(\Omega^\bullet_{\tilde X|_s}(\log \tilde D|_s), 
\tilde\nabla_s)$
is given by the complex $A^{\bullet}_s$. Moreover, if we assume that
$a_i$, $a_h$ and $\sum_{i\in I_L}a_i$ for all $L\in \cal(S\times \arrA)$
are not in $\Z$ one has that this cohomology is equal to the 
cohomology of the local system defined by $ker(\nabla_s)$, 
see [D1, II.6]. We know, from [ESV], that the cohomology of this 
local system is concentrated in degree $n$, see 
Theorem \ref{theorem vanishing local systems}.
We have the following theorem:

\begin{Theorem}
\label{teo:esv-ext NNC}
Let $\XX , S$ and $\pi : \XX \longrightarrow S$ as before. Then,
if
$a_i \notin \Z$  for 
$i \in \{0,...,m,h\}$ with $a_0 = -\sum_{i=1}^m a_i - a_h$ and if
$\sum_{i \in I_L} a_i \not\in \Z$ for every $L \in \calL$, then we
have that ${\bf R}^p(\pi\circ\rho)_*\Omega^\bullet_{\tilde X/S}
(\log \tilde D)=0 \hbox{ \ for \ } p \not= n$ and $A^p=
(\pi\circ\rho)_*\Omega^p_{\tilde X/S}
(\log \tilde D)$.
\end{Theorem} 

\prf
For $\Omega^q_{\tilde X}(\log \tilde D)$ we define the function 
$h^p$ on $S$ as 
\[
h^p(s,\Omega^q_{\tilde X/S}(\log \tilde D))=
\dim H^p(\tilde X|_s,\Omega^q_{\tilde X/S}(\log \tilde D)|_s)
\]
where  $\Omega^q_{\tilde X/S}(\log \tilde D)|_s$
is just the restriction to the fiber $\pi^{-1}(s) = s\times \pn$.
As from Remark \ref{remark familie's bad loci} we have that 
$\Omega^q_{\tilde X/S}(\log \tilde D)|_s = 
\Omega^q_{\tilde X|_s}(\log \tilde D|_s)$,
from Lemma \ref{lemmaVanishing} we have that
$H^p(X,\Omega^q_{\tilde X|_s}
(\log \tilde D|_s)) =0$, for $p>0$.
This implies that for $p>0$ $h^p$ is the constant function with value zero.  
As $S$ is the complement of the discriminant, from Theorem
\ref{teo:nbc} we have that $h^0$ is constant equal 
to $|nbc(\arrA)|$.
We conclude that $h^p$ is constant for every $p$.

As $\pi\circ\rho$ is projective and $h^p$ is constant, 
from [Ha, III, 12.9], we have that for every $s\in S$ the natural
map
\begin{equation}
\label{base change for cohomology NNC}
{\bf R}^p(\pi\circ\rho)_*\Omega^q_{\tilde X/S}(\log \tilde D)
\otimes k(s)\flecha 
    H^p(\tilde X|_s,\Omega^q_{\tilde X/S}(\log \tilde D)|_s)
\end{equation}
where $k(s)$ is the residue field over $s$, is an isomorphism for every $p$.
By Lemma~\ref{lemmaVanishing} we have 
$H^p(\tilde X\mid_s,\Omega^q_{\tilde X/S}(\log \tilde D)|_s)=0$ 
for $p>0$
which from (\ref{base change for cohomology NNC})
implies that
\begin{equation}
\label{eq: vanishing direct images NNC}
{\bf R}^p(\pi\circ\rho)_*\Omega^q_{\tilde X/S}(\log \tilde D) = 0
\end{equation}
for $p>0$.

For $p=0$ we have that 
$H^0(\tilde X\mid_s,\Omega^q_{\tilde X/S}(\log \tilde D)|_s)$
has constant dimension equal to $dim A^p_s$ for every $s\in S$. Thus 
by base change and [Ha, III, 12.9] the inclusion
\[
A^p \hookrightarrow \pi_*\Omega^p_{\smXS}(\log D)
\]
is an isomorphism.

From the spectral sequence
\[
E_1^{p,q}={\bf R}^p(\pi\circ\rho)_*\Omega^q_{\tilde X/S}(\log \tilde D)
\Longrightarrow
{\bf R}^{p+q}(\pi\circ\rho)_*\Omega^\bullet_{\tilde X/S}(\log \tilde D)
\]
and (\ref{eq: vanishing direct images NNC}) we have that
\begin{equation}
\label{eq: coherence NNC}
\begin{array}{c}
{\bf R}^p(\pi\circ\rho)_*\Omega^\bullet_{\tilde X/S}(\log \tilde D) =\\
\\
\underline{ker((\pi\circ\rho)_*\Omega^p_{\tilde X/S}(\log \tilde D) 
   \stackrel{d+\tilde \omega}{\flecha} 
       (\pi\circ\rho)_*\Omega^{p+1}_{\tilde X/S}(\log \tilde D))}\\
im((\pi\circ\rho)_*\Omega^{p-1}_{\tilde X/S}(\log \tilde D)
   \stackrel{d+\tilde \omega}{\flecha}
      (\pi\circ\rho)_*\Omega^p_{\tilde X/S}(\log \tilde D)).
\end{array}
\end{equation}
As $d_{rel} = 0$ on 
$(\pi\circ\rho)_*\Omega^\bullet_{\tilde X/S}(\log \tilde D)$, the 
differential $d_{rel}+\omega = \omega$ and
${\bf R}^p(\pi\circ\rho)_*\Omega^\bullet_{\tilde X/S}(\log \tilde D)$ is
coherent. 

If we denote the maximal ideal of $\hol_{S,s}$ as $\calM_s$, 
from (\ref{eq: coherence NNC}), one can show that one has the base
change formula 
\begin{equation}
\label{projection formula NNC}
{\bf R}^p(\pi\circ\rho)_*((\pi\circ\rho)^*\calM_s\otimes
\Omega^\bullet_{\tilde X/S}(\log \tilde D))
=\calM_s\otimes{\bf R}^p(\pi\circ\rho)_*
\Omega^\bullet_{\tilde X/S}(\log \tilde D).
\end{equation}

We have the following exact sequence
\begin{equation}
\label{restriction sequence NNC}
\begin{array}{c}
0\flecha (\pi\circ\rho)^*\calM_s \otimes 
\Omega^\bullet_{\tilde X/S}(\log \tilde D)\flecha \\
\\
\Omega^\bullet_{\tilde X/S}(\log \tilde D) 
\flecha \Omega^\bullet_{\tilde X/S}(\log \tilde D)|_s \flecha 0.
\end{array}
\end{equation}

For the long exact sequence of cohomology, from 
Theorem~\ref{theorem vanishing local systems}, we have that
for $p\neq n$
\[
{\bf R}^p(\pi\circ\rho)_*(\Omega^\bullet_{\tilde X/S}(\log \tilde D)|_s)=0.
\]
Applying projection formula (\ref{projection formula NNC}) 
gives the following surjection
\[
\calM_s\otimes {\bf R}^p(\pi\circ\rho)_*
\Omega^\bullet_{\tilde X/S}(\log \tilde D)
\flecha {\bf R}^p(\pi\circ\rho)_*
\Omega^\bullet_{\tilde X/S}(\log \tilde D) \flecha 0
\]
for $p \neq n$. By Nakayama's lemma gives
\[
{\bf R}^p(\pi\circ\rho)_*\Omega^\bullet_{\tilde X/S}(\log \tilde D) =0 
\]
for $p \neq n$.

For $p=n$, by (\ref{eq: coherence NNC}) we have that 
${\bf R}^n(\pi\circ\rho)_*\Omega^\bullet_{\tilde X/S}(\log \tilde D)=
A^n/\omega A^{n-1}$ is a coherent $\hol_S$-module.
\qed

\begin{Remark} The sheaf of
$\hol_S$-modules 
${\bf R}^n(\pi\circ\rho)_*\Omega^\bullet_{\tilde X/S}(\log \tilde D)$ 
is free over $\hol_S$.
\end{Remark}

Let $A^{\prime^p}\subset A^p$ be the subalgebra generated by 
$\{\bigwedge_{j=1}^p \frac{dx_{i_j}}{x_{i_j}} ~| \hbox{\ }i_j \in
\{1,...,m\} \hbox{\ and \ } x_{i_j} < x_{i_k} \hbox{\ if\ } j<k
\}$.
On every fiber we have the arrangement $\arrA \cup H_s$, where
$H_s$ is 
defined by $x_s:=1+\sum_{i=1}^nl_ix_i=0$ with 
$s=(1,l_1,...,l_n) \in S \subset \P^{n^{\vee}}$. As for every 
$l \in \calL (\arrA \cup H_s)$, where $\calL$ is as in 
definition \ref{def. bad}, we have that $l \not\subset H_s$
which from Proposition \ref{prop:rel} implies that
\begin{equation}
\label{decomposition local sections NNC}
A^p_s = A^{\prime'^p}_s\oplus \frac{dx_s}{x_s}\wedge A^{\prime^{p-1}}_s.
\end{equation}
From claim \ref{generator of A^p} we have that $A^p_s$ generates
$H^p(\tilde X|_s,\Omega^q_{\tilde X/S}(\log \tilde D)|_s)$.
As a consequence of Theorem \ref{teo:esv-ext NNC}, from 
(\ref{base change for cohomology NNC}) we can extend  
the decomposition (\ref{decomposition local sections NNC}) to
global sections as
\begin{equation}
\label{decomposition global sections NNC}
A^p = A^{\prime^p}\oplus \frac{dx_s}{x_s}\wedge A^{\prime^{p-1}}.
\end{equation}
Applying the Euler characteristic to the sequence 
\[
0\flecha A^0\flecha A^1\flecha\dots\flecha A^n\flecha 
{\bf R}^n(\pi\circ\rho)_*\Omega^\bullet_{\tilde X/S}(\log \tilde D)\flecha 0
\]
we have proved the following corollary.

\begin{Corollary}
\label{cor:basis local systems NNC} 
Under the same assumptions as in Theorem \ref{teo:esv-ext NNC} we
have that $A^{\prime^n}$ generates 
${\bf R}^n(\pi\circ\rho)_*\Omega^\bullet_{\tilde X/S}(\log \tilde D)$ 
where $A^{\prime^n} \subset A^n$ is generated by 
$\{\wedge_{j=1}^p\frac{dx_{i_j}}{x_{i_j}} \mid i_j \in \{1,...,m\}\}$
for which from Theorem \ref{teo:nbc} the set $nbc(\arrA)$ forms a basis.
\end{Corollary}

\begin{Remark}
When the arrangement $\arrA$ has normal crossings the sheaf 
of $\hol_S$-modules ${\bf R}^n\pi_*\Omega^\bullet_{ \smXS}(\log { D})$ 
is free of rank ${m \choose n}$ over $\hol_S$ with basis 
$\{\bigwedge^n_{j=1}\frac{dx_{i_j}}{x_{i_j}}~|~ i_j \in 
\{1,\ldots,n+r\}~ and ~i_j<i_k ~when ~j<k\}$.
\end{Remark}

As under the hypothesis of Theorem \ref{teo:esv-ext NNC}, 
the cohomology of the complex 
$\Omega^\bullet_{\tilde X/S}{\log(\tilde D)}$ 
is concentrated in degree $n$, the Gau\ss-Manin connection 
is given as
\begin{equation}
\label{eq:gm NNC}
{\overline\nabla} : {\bf R}^n(\pi\circ\rho)_*\Omega^\bullet_{ \tilde X/S}
(\log { \tilde D}) \longrightarrow \Omega^1_S \otimes {\bf
R}^n(\pi\circ\rho)_*\Omega^\bullet_{ \tilde X/S} 
(\log  \tilde D).
\end{equation}  

\begin{Remark}
The order on the set of hyperplanes induces an order ``$<$'' 
on the elements of $A^1$. 
This order induces an order on the basis of Corollary 
\ref{cor:basis local systems NNC} for 
${\bf R}^n(\pi\circ\rho)_*\Omega^\bullet_{ \tilde X/S}(\log  \tilde D)$, 
where we say that 
$\bigwedge_{k=1}^n\frac{dx_{i_k}}{x_{i_k}}\leq
\bigwedge_{k=1}^n\frac{dx_{j_k}}{x_{j_k}}$ when there is $k\in\{1,\dots,n\}$
such that $\frac{dx_{i_k}}{x_{i_k}}\leq\frac{dx_{j_k}}{x_{j_k}}$
and $\frac{dx_{i_l}}{x_{i_l}}=\frac{dx_{j_l}}{x_{j_l}}$ for $l<k$.
\end{Remark}

The procedure to write matrix of the Gau\ss-Manin connection with
respect to the basis given in corollary \ref{cor:basis local
systems NNC} is as follows:
We take, as before, affine coordinates for the complement of
$z_0=0$ in $\pn$ as $x_i=z_i/z_0$. We do the same for the complement of 
$h_0=0$ in $\P^{n^{\vee}}$
by taking $l_i=h_i/h_0$. We extend the relative differential form 
$\omega$ to a global form $\Omega$, as in Section 
\ref{section GM-connection}. 
We have that in affine coordinates 
$\Omega = \sum_{i=1}^n a_i\frac{dx_i}{x_i} + a_h\frac{dx_l}{x_l}$ where 
$x_l=1+l_1x_1+...+l_nx_n$ and where the differential is the
absolute one.
The procedure is the standard one, we take an element of the 
basis given in \ref{cor:basis local systems NNC}, we apply to it
the connection and write its image again in terms of this basis. 
With help of the basis of relations given in Proposition \ref{prop:rel} 
one can write the image canonically back in terms of the basis.
We have basically to cases. The first one is when, by taking 
an element of the basis \ref{cor:basis local systems NNC}, the
hyperplanes involved with this elements are given by the set of
affine coordinates chosen. Under the order induced to the basis
this is the first element of our basis.  For the rest of the 
elements we suggest to make a change of basis for the affine
coordinates.

As the basis \ref{cor:basis local systems NNC}, for the 
Gau\ss-Manin bundle, depends on the combinatorics of our
arrangement we cannot give and explicit form for the matrix.
Nevertheless, the basis $nbc(\arrA)$ and the basis of
relations in Proposition \ref{prop:rel} are given in a so
precise way that, for any explicit example, with the method
above, one can compute the matrix for the Gau\ss-Manin connection.


\section{Example I}
\label{sec:exI}
In this section we give an example for the method given in
the previous section. We take an arrangement of six lines
in $\P^2$ in general position.  The discriminant in this 
case is Ceva's arrangement.  This configuration has been deeply
studied in [BHH].

Let $\arrA = \cup_{i=0}^3H_i$ be the arrangement in $\P^2$
given as:
\begin{equation}
\label{fiber nnc}
\begin{array}{l}
H_0 := z_0 = 0 \\
H_1 := z_1 = 0 \\
H_2 := z_2 = 0 \\
H_3 := z_3:= z_0+z_1+z_2=0
\end{array}
\end{equation}
\[
\unitlength.5cm
\begin{array}{c}
\begin{picture}(6,10)
\put (2.8,9.2){\line(1,-2){2.9}}
\put (3.2,9.2){\line(-1,-2){2.9}}
\put (0.1,3.7){\line(2,1){5}}
\put (6.1,3.6){\line(-2,1){5.2}}    
\put (1,7) {$H_0$}
\put (4.5,7) {$H_1$}
\put (1.7,4) {$H_2$}
\put (3.5,4) {$H_3$}
\put (2.5,2) {$\small{Fig. 3}$}
\end{picture}

\end{array}
\]
where we can take $z_0, z_1, z_2$ as a local frame for $\P^2$.
In this case, the discriminant is given as $Discr(\arrA) = \{ h_0=0,
h_1=0, h_2=0, h_0-h_1=0, h_0-h_2=0, h_2-h_1=0\}.$  

\begin{equation}
\label{Ceva}
\begin{array}{c}
\unitlength1cm
\begin{picture}(7,5)
\put (1.7,0.7){\line(1,1){2.7}}
\put (2.3,1.9){$\small{H_5}$}
\put (1.4,3.3){\line(2,-1){5}}
\put (2.3,2.9){$\small{H_4}$}
\put (1.8,5.4){\line(1,-2){2.4}}
\put (3.2,1.5) {$\small{H_3}$}
\put (1.5,1){\line(1,0){5}}        
\put (5, 0.5) {$H_2$}
\put (2,0.4){\line(0,1){5}}        
\put (1, 2 ) {$\tiny{H_1}$}
\put (1.7,5.3){\line(1,-1){4.6}}
\put (4.7,2.5) {$H_0$}
\put (3, 0){$Fig.~4$} 
\end{picture}
\end{array}
\end{equation}

Let $X=S\times\P^2\setminus
\{\Delta:=h_0z_0+h_1z_1+h_2z_2=0\}\cup\{S\times\arrA\}$ and $\pi: 
X\flecha S$ be a family 
of arrangements parametrized by
$S={\P^2}^\vee\setminus Discr(\arrA)$. We denote the divisor 
$S\times\arrA\cup\{\Delta\cap S\times\P^2 \}$ by $D$. 

We fix $H_0$ as the hyperplane at infinity of the arrangement 
(\ref{fiber nnc}). 
Let $\omega \in 
H^0(S\times\P^2, \Omega^1_{S\times\P^2}(\log D))$ be given as
$\omega = \sum_{i=1}^3 a_i\frac{d_{rel}x_i}{x_i}+a_h\frac{d_{rel}x_l}{x_l}$
where $x_i=z_i/z_0$ and $l_i=h_i/h_0$ and $x_l=l_1x_1+l_2x_2+1$,
and where the 
differential is taken as the relative differential along $S$, and 
where $\frac{dx_i}{x_i}$ is taken as in Remark \ref{notation ueber die
Differentiale Formen}. We assume that $a_i \not\in \Z$ for $i \in \{~
0,\dots,3, h~\}$, and that
$ \sum_{i=0}^3 a_i+a_h = 0$.

The operator $\nabla = d_{rel} + \omega$ defines the complex
\[ 
0\flecha \hol_X \flecha \Omega_{\smXSdos}^1(\log D) \flecha 
\Omega_{\smXSdos}^2(\log D) \flecha  0.
\]
Let $V$ be the relative local system defined as the sheaf of flat 
sections of $\nabla$.
Under the above assumptions on $\omega$, from \ref{teo:esv-ext NNC} 
we have that the cohomology
\begin{equation}
H^i(X/S,V) = {\bf R}^i\pi_*\Omega_{\smXSdos}^{\bullet}(\log D) =0
\end{equation}
for  $i = 0,1$.

For $i=2$ we have 
\begin{equation}
\label{eq:baseex1}
\begin{array}{c}
H^2(X/S,V)={\bf R}^2\pi_*\Omega_{\smXSdos}^\cdot (\log D)=\\
\\
\hol_S\frac{dx_1\wedge dx_2}{x_1x_2}\oplus 
\hol_S\frac{dx_1\wedge dx_3}{x_1x_3}\oplus
\hol_S\frac{dx_2\wedge dx_3}{x_2x_3}
\end{array}
\end{equation}

We can now extend $\omega$ to
\[
\Omega = \sum_{i=0}^3 a_i\frac{dz_i}{z_i}+a_h\frac{dz_h}{z_h}
\]
where the differential is no longer the relative differential. We have 
the operator ${\overline\nabla} = d + \Omega$ which, when using affine
coordinates namely on the complement of $z_0=0$ and $h_0=0$, takes the
form 
\[
{\overline\nabla} = d +\sum_{i=1}^3 a_i\frac{dx_i}{x_i}+a_h\frac{dx_l}{x_l}.
\]
The Gau\ss-Manin connection is then obtained with help of the canonical
filtration applied to the de~Rham complex 
$(\Omega_{\smX}^\bullet(\log D),{\overline\nabla})$ and is given as
\[
{\overline\nabla} : H^2(X/S,V) \flecha H^2(X/S,V)\otimes 
\Omega^1_S(\log(Discr(\arrA))).
\]

With respect to the basis (\ref{eq:baseex1}) one can represent
this connection by the matrix
\begin{equation}
\label{GM matrix ex. 1}
\small{
\left(
\begin{array}{ccc}
    \begin{array}{c}
        -a_1[\frac{dh_1}{h_1}-\frac{dh_0}{h_0}]\\
        ~~~-a_2[\frac{dh_2}{h_2}-\frac{dh_0}{h_0}] 
    \end{array}       
            & -a_2[\frac{dh_2}{h_2}-\frac{d(h_0-h_2)}{h_0-h_2}] 
                   & a_1[\frac{dh_1}{h_1}-\frac{d(h_0-h_1)}{h_0-h_1}]   \\
&&\\
&&\\
-a_3[\frac{dh_2}{h_2}-\frac{dh_0}{h_0}] 
       & \begin{array}{c}
             -a_1[\frac{d(h_1-h_2)}{h_1-h_2}-\frac{d(h_0-h_2)}{h_0-h_2}]\\
              ~~~~~~-a_3[\frac{dh_2}{h_2}-\frac{d(h_0-h_2)}{h_0-h_2}]  
         \end{array}
       & -a_1[\frac{d(h_1-h_2)}{h_1-h_2}-\frac{d(h_0-h_1)}{h_0-h_1}]\\
&&\\
&&\\
a_3[\frac{dh_1}{h_1}-\frac{dh_0}{h_0}]
       & -a_2[\frac{d(h_1-h_2)}{h_1-h_2}-\frac{d(h_0-h_2)}{h_0-h_2}] 
       & \begin{array}{c}
             -a_2[\frac{d(h_1-h_2)}{h_1-h_2}-\frac{d(h_0-h_1)}{h_0-h_1}]\\
              ~~~~~~~-a_3[\frac{dh_1}{h_1}-\frac{d(h_0-h_1)}{h_0-h_1}] 
         \end{array} 
\end{array}
\right).}
\end{equation}

We would now like to calculate the monodromy of the Gau\ss-Manin
connection along different elements of the fundamental group
of $S$.

For $H_i \in Discr(\arrA)$ we have the residue map along $H_i$
\[
Res_{H_i}({\overline\nabla}) : H^2(X/S,V) \flecha H^2(X/S,V)\otimes 
\Omega^1_S(\log Discr(\arrA)) \flecha H^2(X/S,V)\otimes \hol_{H_i}
\]
and defined in the usual way, see [D1, II.3.7].

Fix a base point $p \in S$ and let 
$\gamma_i \in \pi(S,p)$ be a loop around $H_i \in Discr(\arrA)$
with base point $p$. 

Let 
\begin{equation}
\label{equation for monodromy representation}
T_i= \exp (-2\pi i\cdot Res_{H_i}(\overline\nabla)).
\end{equation}
If we suppose that $a_i+a_j \not\in\Z\setminus \{0\}$ for 
$1\leq i<j\leq 3$ then, by [D1, II.5.6], the monodromy transformation 
$M_p$ when we go around $H_i$ along $\gamma_i$ is given by a conjugacy
class of $T_i$.


The monodromy around $H_1 : h_1 = 0$ is given as follows.

Let $A_{H_1}$ be the residue matrix of the connection along $H_1$.
From (\ref{GM matrix ex. 1}) we have that
\[
A_{H_1} = Res_{H_1}(\overline\nabla) =
\left(
\begin{array}{ccc}
-a_1 & 0 & a_1 \\
  0  & 0 &  0  \\
 a_3 & 0 & -a_3
\end{array}
\right).
\]

We have that for $n \geq 1$
\[
A_{H_1}^n = (-a_1-a_3)^{n-1}A_{H_1}.
\]
One can see that $(-a_1-a_3)$ is the trace of the matrix $A_{H_1}$ 
which is an eigenvalue.
We have that $A_{H_1}^n = tr(A_{H_1})^{n-1}A_{H_1}$
where $tr$ is the trace of the matrix.
If $a_1+a_3 \not\in \Z\setminus\{0\}$ from 
(\ref{equation for monodromy representation}) we have then 
that the monodromy transform is given by a conjugacy class of
\[ 
T_1 = I+ (\exp(-2\pi i\cdot (-a_1-a_3)) -1)\cdot
\left(
\begin{array}{ccc}
\frac{-a_1}{(-a_1-a_3)} & 0 & \frac{a_1}{(-a_1-a_3)} \\
              0         & 0 &                 0  \\
\frac{a_3}{(-a_1-a_3)}  & 0 & \frac{-a_3}{(-a_1-a_3)}
\end{array}
\right).
\]
This is
\[  
M_p(\gamma_1)= I + (\exp(-2\pi i \cdot tr(A_{H_1})) - 1) 
                      tr(A_{H_1})^{-1}\cdot \alpha A_{H_1}\alpha^{-1}.
\]

For $H_2 : h_2  = 0$ we have the following.

Let $A_{H_2}$ be the residue matrix of the connection 
along $H_2$.
From (\ref{GM matrix ex. 1}) we have that
\[
A_{H_2} = Res_{H_2}(\overline\nabla) =
\left(
\begin{array}{ccc}
-a_2 & -a_2 &  0 \\
-a_3 & -a_3 &  0  \\
  0  &  0   &  0
\end{array}
\right).
\]

We have that for $n \geq 1$
\[
A_{H_2}^n = (-a_2-a_3)^{n-1}A_{H_2}.
\]
Again we have that $A_{H_2}^n = tr(A_{H_2})^{n-1}A_{H_2}$
where $tr$ is the trace and an eigenvalue of the matrix $A_{H_2}$.
If $a_2+a_3 \not\in \Z\setminus\{0\}$, from 
(\ref{equation for monodromy representation}) we have then 
that the monodromy transform is given by a conjugacy class of
\[  
T_2= I + (\exp(-2\pi i\cdot tr(A_{H_2})) - 1)
                  tr(A_{H_2})^{-1}\cdot A_{H_2}.
\]

We calculate now the monodromy around $H_0 : h_0 = 0$.

As before let $A_{H_0}$ be the residue matrix of the 
connection along $H_0$.
From (\ref{GM matrix ex. 1}) we have that
\[
A_{H_0} = Res_{H_0}(\overline\nabla) =
\left(
\begin{array}{ccc}
 a_1+a_2  &  0  &  0  \\
    a_3   &  0  &  0  \\
   -a_3   &  0  &  0
\end{array}
\right).
\]

We have that for $n \geq 1$
\[
A_{H_0}^n = (a_1+a_2)^{n-1}A_{H_0}.
\]
Again we see that $A_{H_0}^n = tr(A_{H_0})^{n-1}A_{H_0}$
where $tr$ is the trace and an eigenvalue of the matrix $A_{H_0}$.
If $a_1+a_2 \not\in \Z\setminus\{0\}$,
from (\ref{equation for monodromy representation}) we have then 
that the monodromy transform is given as a conjugacy class of
\[  
T_0= I + (\exp(-2\pi i\cdot tr(A_{H_0})) - 1)
                  tr(A_{H_0})^{-1}\cdot A_{H_0}.
\]

The monodromy around $H_3 : h_0-h_1=0$ is given as follows.

Let $A_{H_3}$ be the residue matrix of the connection along $H_3$.
From (\ref{GM matrix ex. 1}) we have that
\[
A_{H_3} = Res_{H_3}(\overline\nabla) =
\left(
\begin{array}{ccc}
  0  &  0  & -a_1 \\
  0  &  0  &  a_1  \\
  0  &  0  & a_2+a_3
\end{array}
\right).
\]

For $n \geq 1$ we have that
\[
A_{H_3}^n = (a_2+a_3)^{n-1}A_{H_3}.
\]
One can see that $(a_2+a_3)$ is the trace and an eigenvalue of the 
matrix $A_{H_3}$.
We have again that $A_{H_3}^n = tr(A_{H_3})^{n-1}A_{H_3}$
where $tr$ is the trace of the matrix.
If $a_1+a_3 \not\in \Z\setminus\{0\}$, 
from (\ref{equation for monodromy representation}) we have then 
that the monodromy transform is given as a conjugacy class of
\[  
T_3= I + (\exp(-2\pi i\cdot tr(A_{H_3})) - 1)
                 tr(A_{H_3})^{-1}\cdot A_{H_3}.
\]

Around $H_4 : h_0-h_2 = 0$ the monodromy is given as follows.

Let $A_{H_4}$ be the residue matrix of the connection along $H_4$.
From (\ref{GM matrix ex. 1}) we have that
\[
A_{H_4} = Res_{H_4}(\overline\nabla) =
\left(
\begin{array}{ccc}
  0  &    a_2   &  0  \\
  0  &  a_1+a_3 &  0  \\
  0  &    a_2   &  0
\end{array}
\right).
\]

We have that for $n \geq 1$
\[
A_{H_4}^n = (a_1+a_3)^{n-1}A_{H_4}.
\]
One can see again that $(a_1+a_3)$ is the trace and an eigenvalue 
of the matrix $A_{H_4}$.
We have that $A_{H_4}^n = tr(A_{H_4})^{n-1}A_{H_4}$
where $tr$ is the trace of the matrix.
If $a_1+a_3 \not\in \Z\setminus\{0\}$, 
from (\ref{equation for monodromy representation}) we have then 
that the monodromy transform is given as a conjugacy class of
\[  
T_4= I + (\exp(-2\pi i\cdot tr(A_{H_4})) - 1)
                 tr(A_{H_4})^{-1}\cdot A_{H_4}.
\]

The monodromy around $H_5 : h_2-h_1 = 0$ is given as follows.

Let $A_{H_5}$ be the residue matrix of the connection along $H_5$.
From (\ref{GM matrix ex. 1}) we have that
\[
A_{H_5} = Res_{H_5}(\overline\nabla) =
\left(
\begin{array}{ccc}
  0  &    0    &  0   \\
  0  &  -a_1   & -a_1  \\
  0  &  -a_2   & -a_2
\end{array}
\right).
\]

We have that for $n \geq 1$
\[
A_{H_5}^n = (-a_1-a_2)^{n-1}A_{H_5}.
\]
One can see that $(-a_1-a_2)$ is the trace and an eigenvalue 
of the matrix $A_{H_5}$.
We have that $A_{H_5}^n = tr(A_{H_5})^{n-1}A_{H_5}$
where $tr$ is the trace of the matrix.
If $a_1+a_2 \not\in \Z\setminus\{0\}$, 
from (\ref{equation for monodromy representation}) we have then 
that the monodromy transform is given as a conjugacy class of
\[  
T_5= I + (\exp(-2\pi i\cdot tr(A_{H_5})) - 1)
                 tr(A_{H_5})^{-1}\cdot A_{H_5}.
\]


\section{Ceva's Configuration}
\label{Ceva's Configuration}

Let $\arrA = \cup_{i=0}^5H_i$ be an arrangement in $\P^2$
given as:
\begin{equation}
\label{projective Ceva discriminant}
\begin{array}{l}
H_0 := z_0 = 0 \\
H_1 := z_1 = 0 \\
H_2 := z_1 = 0 \\
H_3 := z_3:= z_0-z_1=0\\
H_4 := z_4:= z_0-z_2=0\\
H_5 := z_5:= z_1-z_2=0
\end{array}
\end{equation}
\[
\begin{array}{c}
\unitlength1cm
\begin{picture}(7,5)
\put (1.7,0.7){\line(1,1){2.7}}
\put (2.3,1.9){$\small{H_5}$}
\put (1.4,3.3){\line(2,-1){5}}
\put (2.3,2.9){$\small{H_4}$}
\put (1.8,5.4){\line(1,-2){2.4}}
\put (3.2,1.5) {$\small{H_3}$}
\put (1.5,1){\line(1,0){5}}        
\put (5, 0.5) {$H_2$}
\put (2,0.4){\line(0,1){5}}        
\put (1, 2 ) {$\tiny{H_1}$}
\put (1.7,5.3){\line(1,-1){4.6}}
\put (4.7,2.5) {$H_0$}
\put (3, 0){$Fig.~5$} 
\end{picture}
\end{array}
\]

The discriminant is $Discr(\arrA) = \{ h_0= 0, h_1=0, h_2=0,
h_0+h_1=0, h_1+h_2=0, h_0+h_2=0, h_0+h_1+h_2=0 \}$.

Let $X=S\times\P^2\setminus
\{\Delta:=h_0z_0+h_1z_1+h_2z_2=0\}\cup\{S\times\arrA\}$ and $\pi:X \flecha S$ 
be a family of arrangements parametrized by
$S={\P^2}^\vee\setminus Discr(\arrA)$. We denote the divisor 
$(S\times\arrA)\cup\{\Delta\cap (S\times\P^2) \}$ as $D$. 

Let $\rho : \tilde X \flecha S\times\P^2$ be the blow up along
the elements of $\calL(D)$ as in (\ref{family normalization}),
see remark \ref{remark familie's bad loci}. Let $\tilde D = \rho^*(D)$.

Let $H_0$ be the hyperplane at infinity of the projective arrangement
(\ref{projective Ceva discriminant}). The affine arrangement
on the complement of $H_0$ is given then by the equations 
$\{x_1=0, x_2=0, x_3:=x_1-1=0, x_4:= x_2-1=0, x_5:=x_1-x_2=0\}$.
Let $U = S\times\P^2 \setminus D$. Let $\omega \in 
H^0(U, \Omega^1_U)$ be given as
$\omega = \sum_{i=1}^5 a_i\frac{d_{rel}x_i}{x_i}+
a_h\frac{d_{rel}x_l}{x_l}$
where $x_i=z_i/z_0$ and $l_i=h_i/h_0$ and $x_l=l_1x_1+l_2x_2+1$.
and where $\frac{d_{rel}x_i}{x_i}$ is taken as in remark \ref{notation
ueber die Differentiale Formen} and where the differential is taken as
the relative differential along $S$. We assume that 
$ \sum_{i=0}^5 a_i+a_h = 0$, $a_i \not\in \Z$ for $i \in 
\{~0,\dots,5,h~\}$ and 
$\sum_{i\in I_L}a_i \not\in \Z$ for $L\in \calL(\arrA)$.
Let $\tilde \omega \in H^0(\tilde X, \Omega^1_{\tilde X}(\log \tilde D))$
be given as $\tilde \omega = \rho^* \omega$.

Again as in section \ref{sec:exI} we have that the operator 
$\bar\nabla = d_{rel} + \tilde\omega$ defines the complex
\[ 
0\flecha \hol_{\tilde X} \flecha \Omega_{\tilde X}^1(\log \tilde D) 
\flecha \Omega_{\tilde X/S}^2(\log \tilde D) \flecha  0.
\]
From Theorem \ref{teo:esv-ext NNC} and [D1, II.6] we have that 
the cohomology of the local system $V$ obtained as the flat sections 
of $\bar\nabla$ is
\begin{equation}
H^i(\tilde X/S,V) = 
{\bf R}^i\pi_*\Omega_{\tilde X/S}^{\bullet}(\log \tilde D) = 0
\end{equation}
for  $i = 0,1$.

For $i=2$ we have from corollary \ref{cor:basis local systems NNC}
\begin{equation}
\label{basis for cohomology for projective ceva}
\begin{array}{c}
H^2(\tilde X/S,V)=
       {\bf R}^2\pi_*\Omega_{\tilde X/S}^\cdot (\log \tilde D)=\\
\\
\hol_S\frac{dx_1\wedge dx_2}{x_1x_2}\oplus 
\hol_S\frac{dx_1\wedge dx_4}{x_1x_4}\oplus
\hol_S\frac{dx_1\wedge dx_5}{x_1x_5}\oplus
\hol_S\frac{dx_2\wedge dx_3}{x_2x_3}\oplus
\hol_S\frac{dx_3\wedge dx_4}{x_3x_4}\oplus
\hol_S\frac{dx_3\wedge dx_5}{x_3x_5}
\end{array}
\end{equation}

We have the basis of relations of elements of $A^2$ given in 
(\ref{relations ceva projective}).

Let $\Omega \in H^0(U,\Omega^1_U)$ be given by
\[
\Omega = \sum_{i=1}^5 a_i\frac{dx_i}{x_i}+a_h\frac{dx_l}{x_l}
\]
where the differential is no longer the relative differential
and where $x_i=z_i/z_0$, $l_i=h_i/h_0$, $x_l = 1+l_1x_1+l_2x_2$ and
$dx_i/x_i$ is taken as in remark \ref{notation ueber die Differentiale
Formen}. We extend  $\tilde\omega$ to $\tilde X$ to an 
element $\tilde \Omega 
\in H^0(\tilde X, \Omega^1_{\tilde X}(\log \tilde D))$ as 
$\tilde \Omega = \rho^*\Omega$. We have 
the operator ${\overline\nabla} = d + \tilde \Omega$.  
We have the Gau\ss-Manin connection
\[
\overline\nabla : H^2(\tilde X/S,V) \flecha 
H^2(\tilde X/S,V)\otimes \hol_S.
\]

To write matrix of the Gau\ss-Manin connection with respect to the
basis (\ref{basis for cohomology for projective ceva}) we
apply the relations \ref{relations ceva projective}.
The matrix is given as follows.

The first column is
\begin{equation}
\left(
\begin{array}{c}
(-a_1-a_5)[\frac{dh_1}{h_1}-\frac{dh_0}{h_0}]
         -a_2[\frac{dh_2}{h_2} -\frac{dh_0}{h_0}]\\
\\
-a_4[\frac{dh_2}{h_2}-\frac{dh_0}{h_0}]\\
\\
a_5[\frac{dh_1}{h_1} - \frac{dh_2}{h_2}] \\
\\
a_3[\frac{dh_1}{h_1}-\frac{dh_0}{h_0}]\\
\\
          0             \\
\\
          0              
\end{array}
\right).
\end{equation}

The second column is
\[
\left(
\begin{array}{c}
-a_2[\frac{dh_2}{h_2}-\frac{d(h_0+h_2)}{(h_0+h_2)}]
\\ \\
  -a_1[\frac{dh_1}{h_1}-\frac{d(h_0+h_2)}{(h_0+h_2)}]
         -a_4[\frac{dh_2}{h_2}-\frac{d(h_0+h_2)}{(h_0+h_2)}] 
\\ \\
-a_5[\frac{dh_2}{h_2}-\frac{d(h_0+h_2)}{(h_0+h_2)}]
\\ \\
0
\\ \\
(-a_3-a_5)[\frac{dh_1}{h_1}-\frac{d(h_0+h_2)}{(h_0+h_2)}]
\\ \\
a_5[\frac{dh_1}{h_1}-\frac{d(h_0+h_2)}{(h_0+h_2)}]
\end{array}
\right). 
\]

The third column is

\[
\left(
\begin{array}{c}
a_2[\frac{d(h_1+h_2)}{h_1+h_2}-\frac{dh_2}{h_2}]
\\ \\
 -a_4[\frac{dh_2}{h_2}-\frac{dh_0}{h_0} ]
\\ \\
(-a_1-a_2)[\frac{d(h_1+h_2)}{h_1+h_2}-\frac{dh_0}{h_0}]
         -a_5[\frac{dh_2}{h_2}-\frac{dh_0}{h_0}]
\\ \\
0
\\ \\
a_4[\frac{d(h_1+h_2)}{h_1+h_2}-\frac{dh_0}{h_0}]
\\ \\
(-a_3-a_4)[\frac{d(h_1+h_2)}{h_1+h_2}-\frac{dh_0}{h_0}]
\end{array}
\right).
\]

The fourth column is

\[
\left(
\begin{array}{c}
(a_1+a_5)[\frac{dh_1}{h_1}-\frac{d(h_0+h_1)}{h_0+h_1}]
\\ \\
0
\\ \\
(-a_5[\frac{dh_1}{h_1}-\frac{d(h_0+h_1)}{h_0+h_1}]
\\ \\
-a_3[\frac{dh_1}{h_1}-\frac{d(h_0+h_1)}{h_0+h_1}]
      -a_2[\frac{dh_2}{h_2}-\frac{d(h_0+h_1)}{h_0+h_1}]
\\ \\
a_4[\frac{dh_2}{h_2}-\frac{d(h_0+h_1)}{h_0+h_1}]
\\ \\
a_5[\frac{dh_2}{h_2}-\frac{d(h_0+h_1)}{h_0+h_1}]
\end{array}
\right).
\]
\newpage

The fifth column is

\[
\left(
\begin{array}{c}
0
\\ \\
-a_1[\frac{dh_1}{h_1}-\frac{d(h_0+h_1+h_2)}{h_0+h_1+h_2}] 
\\ \\
0
\\ \\
a_2[\frac{dh_2}{h_2} -\frac{d(h_0+h_1+h_2)}{h_0+h_1+h_2}]
\\ \\
(-a_3-a_5)[\frac{dh_1}{h_1}-\frac{d(h_0+h_1+h_2)}{h_0+h_1+h_2}]
      -a_4[\frac{dh_2}{h_2}-\frac{d(h_0+h_1+h_2)}{h_0+h_1+h_2}]
\\ \\
a_5[\frac{dh_1}{h_1}-\frac{dh_2}{h_2}]
\end{array}
\right).
\]

The sixth column is

\[
\left(
\begin{array}{c}
a_2[\frac{d(h_1+h_2)}{h_1+h_2}-\frac{d(h_0+h_1+h_2)}{h_0+h_1+h_2} ]
\\ \\
0
\\ \\
(-a_1-a_2)[\frac{d(h_1+h_2)}{h_1+h_2}-\frac{d(h_0+h_1+h_2)}{h_0+h_1+h_2}] 
\\ \\
a_2[\frac{dh_2}{h_2} -\frac{d(h_0+h_1+h_2)}{h_0+h_1+h_2}]
\\ \\
a_4[\frac{d(h_1+h_2)}{h_1+h_2}-\frac{dh_2}{h_2}]
\\ \\
(-a_3-a_4)[\frac{d(h_1+h_2)}{h_1+h_2}-\frac{d(h_0+h_1+h_2)}{h_0+h_1+h_2}]
           -a_5[\frac{dh_2}{h_2}-\frac{d(h_0+h_1+h_2)}{h_0+h_1+h_2}] 
\end{array}
\right).
\]


\section{Bibliography}

[A] Aomoto, K., {\it Un th\'eor\`eme du type de Matsushima-Murakami 
concernant l'int\'egrale des fonctions multiformes.} J.~Math.~Pures 
Appl.~52, 1-11 (1973)\\
\\ \
[B] Brieskorn, E., {\it Sur les groupes de tresses, S'm. Bourbaki 
401}. Lect. Notes Math., vol.~317, pp.~21-44 Springer 1973.\\
\\ \
[BHH] Barthel, G., Hirzebruch, F., H\"ofer T., {\it 
  Geradenkonfiguarationen und algebraische Fl\"achen}, Vieweg, Aspects
of Mathematic vol.~D4, (1987).\\
\\ \
[Bj] Bj\"orner, A., {\it On the homology of geometric lattices},
Alg.~Universalis, 14 (1982) 107-128.\\
\\ \
[Bo] Bott, R., {\it Homogeneous vector bundles}, Ann. of Math., 66
(1957) 203-248.\\
\\ \
[CH1] Chen, K.-T., {\it Iterated integrals and exponential homomorphisms}, 
 Proc.~London Math.~Soc.~(3), vol..~4 (1954), 502-512.\\
\\ \
[CH2] Chen, K.-T., {\it Integrals of paths, geometric invariants 
and a generalized Baker-Hausdorff formula}, 
Ann.~of Math., vol.~65 (1957) 163-178.\\
\\ \
[D1] Deligne, P., {\it \'Equations diff\'erentielles \`a points 
singuliers r\'eguliers}, Springer Lect. Notes Math. 163 (1970).\\
\\ \
[D2] Deligne, P., {\it Th\'eorie de Hodge, II}, IHES, 
Public.~math.~no.40 (1972), 5-57.\\ 
\\ \
[DM] Deligne, P.~and Mostow G.~D.~, {\it Monodromy of hypergeometric
     functions and non-lattice integral monodromy}, IHES, 
Public.~math.~no.~63 (1986) 5-90.\\
\\ \ 
[EV] Esnault, H., Viehweg, E., {\it Lectures on Vanishings Theorems},
Birkh\"auser DVM Seminar 20 (1992).\\
\\ \
[EV1] Esnault, H., Viehweg, E., {\it Logarithmic De~Rham complexes
and vanishing theorems}, Invent. math. 86, (1986) 161-194.\\     
\\ \ 
[ESV] Esnault, H., Schechtman, V., Viehweg, E., {\it Cohomology
of local systems of the complement of hyperplanes}, Invent.~math. 109
 (1992), 557-561; erratum, ibid.~ 112 (1993)\\
\\ \
[FT] Falk, M., Terao H., {\it $\beta$nbc-bases for cohomology of local 
systems on hyperplane complements }, to appear\\
\\ \ 
\noindent [H] Hain, R., {\it The geometry of the mixed Hodge structure on
the fundamental group}, Algebraic geometry, Bowdoin 1985, Proc.~Symp.~
Pure Math., vol.~46, part 2 (1987), 247-282. \\
\\ \
[Ha] Hartshorne R., {\it Algebraic Geometry}, Springer GTM 52, (1977).\\
\\ \ 
[K] Katz, N., {\it Nilpotent connections and the monodromy theorem:
application of a result of Turritin}, 
Institut des Hautes \'Etudes Scientifique, Public.~math.~no.~31 (1971)
355-411\\
\\ \
[KO] Katz, N., Oda, T., {\it On the differentiation of de Rham cohomology
classes with respect to parameters}, J.~Math.~Kyoto Univ., 8 (1968), 
199-213.\\
\\ \
[LD] Lappo-Danilevski, A., {\it M\'emoires sur la th\'eorie des 
syst\`emes des \'equations diff\'erentielles lin\'eaires}, reprint,
Chelsea, New York, (1953).\\
\\ \
[M] Messing, W.,~{\it On the nilpotence of hypergeometric equations},
J.~Math.~Kyoto Univ.,~vol.~12, (1971), 369-383.\\
\\ \
[O] Oxley, J., {\it Matroid theory}, Oxford Science Public. (1992).\\
\\ \
[OS] Orlik, P., Solomon, L., {\it Combinatorics and topology of 
complements of hayperplanes}, Invent. Math.,56 (1980) 167-189.\\ 
\\ \ 
[OSS] Okonek, Ch., Schneider, M.,~Spindler,~H.,~{\it Vector Bundles on
Complex Projective Spaces}, Birkh\"auser (1980).\\
\\ \
[OT] Orlik, P., Terao, H., {\it Arrangements of hyperplanes}, Springer
GMW 300, (1992).\\
\\ \
[P] Plemelj, J., {\it Problems in the sense of Riemann and Klein},
Wiley, New York, (1964).\\
\\ \
[S] Simpson, T.~{\it Transcendental aspects of the Riemann-Hilbert
Correspondance}, Ill.~J.~Math.,~vol.~34,~num.~2, (1990), 368-391.\\
\\ \
[SV] Schechtman, V., Varchenko, A., {\it Arrangements of
hyperplanes and Lie algebra homology}, Invent.~Math.,106 (1991) 139-194.\\
\\ \

\end{document}